\documentclass[10pt,letterpaper]{article}
\usepackage[top=0.85in,left=2.75in,footskip=0.75in]{geometry}

\usepackage{float}
\usepackage{color}
\usepackage{colortbl}
\usepackage{amsmath,amssymb}

\usepackage{changepage}

\usepackage[utf8x]{inputenc}

\usepackage{textcomp,marvosym}

\usepackage{cite}

\usepackage{nameref}

\usepackage[right]{lineno}

\usepackage{microtype}
\DisableLigatures[f]{encoding = *, family = * }

\usepackage[table]{xcolor}

\usepackage{array}
\usepackage{algorithm}
\usepackage{algpseudocode}
\newcolumntype{+}{!{\vrule width 2pt}}

\newlength\savedwidth

  \raggedright
  \usepackage[aboveskip=1pt,labelfont=bf,labelsep=period,justification=raggedright,singlelinecheck=off]{caption}

  \makeatletter
  \renewcommand{\@biblabel}[1]{\quad#1.}
    \makeatother
    
    \date{}
    
    \usepackage{lastpage,fancyhdr,graphicx}
    \fancyheadoffset[L]{2.25in}
    \fancyfootoffset[L]{2.25in}
\begin{document}
    \vspace*{0.2in}
    
    \begin{flushleft}
    {\Large
      \textbf\newline{Generating global network structures by triad types} 
    }
    \newline
    \\
    Marjan Cugmas\textsuperscript{1},
    Anu{\v{s}}ka Ferligoj\textsuperscript{1},
    Ale{\v{s} \v{Z}}iberna\textsuperscript{1}
    \\
    \bigskip
    \textbf{1} Faculty of Social Sciences, University of Ljubljana, Ljubljana, Slovenia
    \\
    \bigskip
    
    * marjan.cugmas@fdv.uni-lj.si
    
    \end{flushleft}
    \section*{Abstract}
This paper addresses the question of whether it is possible to generate networks with a given global structure (defined by selected blockmodels, i.e., cohesive, core-periphery, hierarchical and transitivity), considering only different types of triads. Two methods are used to generate networks: (i) the method of relocating links; and (ii) the Monte Carlo Multi Chain algorithm implemented in the \texttt{ergm} package implemented in R. Although all types of triads can generate networks with the selected blockmodel types, the selection of only a subset of triads improves the generated networks' blockmodel structure. However, in the case of a hierarchical blockmodel without complete blocks on the diagonal, additional local structures are needed to achieve the desired global structure of generated networks. This shows that blockmodels can emerge based on only local processes that do not take attributes into account.

    \section*{Introduction}
    
 In social network analysis, considerable attention is paid to global network structures, which can be described using a blockmodel \cite{doreian2005generalized}. A blockmodel consists of groups (also called positions) of units and the relationships between those groups. The units are assigned to the same group if they are equivalent according to the pattern of links to the other units. Often, structural equivalence is assumed. Two units are structurally equivalent if they are linked to the same units and by the same others \cite{lorrain1971structural}, implying they share the same social role \cite{luczkovich2003defining}. There are several well-known and studied types of blockmodels, e.g., cohesive, core-periphery, transitivity and hierarchical. 
 
\begin{figure}[!h]
\begin{adjustwidth}{-2.25in}{0in}
\caption{{\bf The collection of all triad types (triad census).} 
The labels consist of three digits: the first digit denotes the number of mutual links ($\leftrightarrow$), the second stands for the number of arcs ($\rightarrow$) while the third denotes the number of missing links between two units. Some types of triads with the same distribution of links are further differentiated (see columns) and labeled with a letter (C stands for cycle, T for transitivity, U for up and D for down).
}
\label{Fig1}
\includegraphics[width=1.4\textwidth]{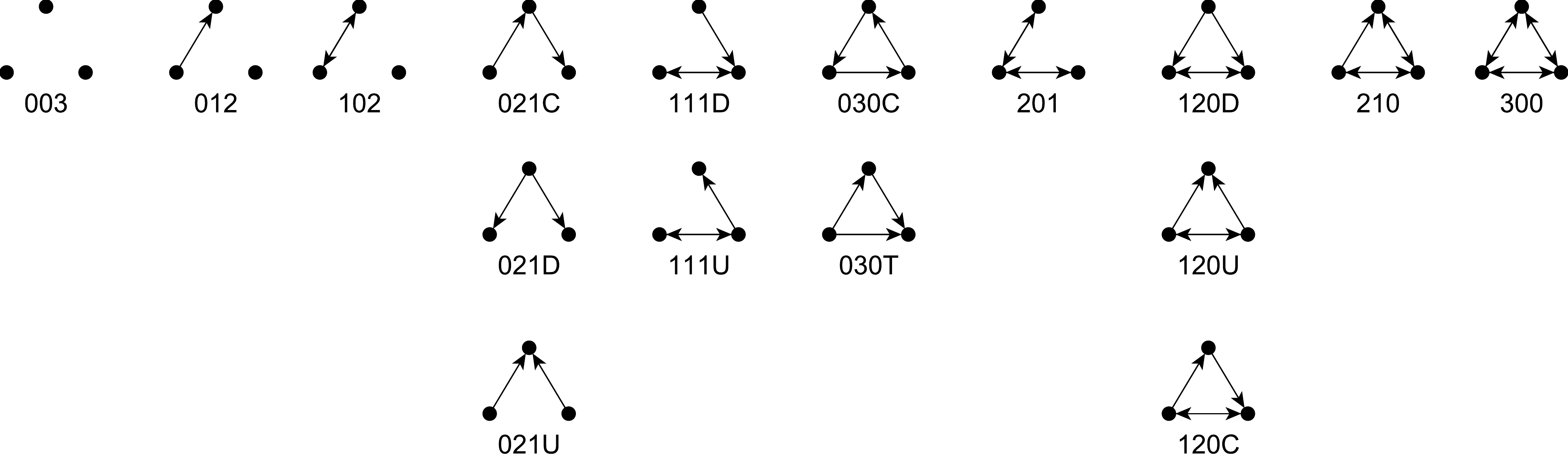}
\end{adjustwidth}
\end{figure}

In order to describe the processes that produce a given global structure, much effort was made to study micro structures in the context of various global structures. For this propose, the triadic census (the collection of all possible networks of size three which are visualised in Fig~\ref{Fig1}), proposed by Davis \cite{davis1967structure}, is often considered \cite{cartwright1956structural, davis1977clustering, davis1967structure, holland1971transitivity, johnsen1985network}.

While the triadic census is well studied in the context of different global network structures, no attention was given to the dependencies between the triadic census and different global network structures, operationalised by the types of blockmodels. \textbf{Therefore, the main objective of this study is to test whether it is possible to generate networks with a given blockmodel structure (seven selected blockmodel types are considered), taking only different types of triads into account.} This is especially important to consider when thinking about the factors that drive a network to a certain global structure in terms of social processes. Understanding such mechanisms might also contribute to the development of statistical tests for assessing the significance of the blockmodel structure obtained by simulations \cite{borgatti2000models} or to predict missing links in partly known networks \cite{clauset2008hierarchical}. In this context, the triad census has already been used to detect the social rules of a participant in a social system \cite{doran2015discovery}, distinguish between network brokerage and network closure \cite{prell2008looking}, predict future links in a social network \cite{juszczyszyn2011link} and study how stature, relationship strength and egocentricity affect user interactions on Facebook \cite{doran2013triads}.
    
The main objective of this study is further elaborated: is it possible to generate networks with a given blockmodel type while considering only allowed or only forbidden types of triads? Therefore, the classification of allowed and forbidden types of triads is determined for each type of blockmodel. Allowed types of triads are those whose frequency is higher than zero in an ideal blockmodel structure. On the other hand, forbidden triad types are those whose frequency equals zero in an ideal blockmodel. The sets of allowed and forbidden triad types are then reduced based on comparisons of the different types of blockmodels for different levels of errors in the network according to the ideal blockmodel being considered. 

The sets of all triad types, the sets of allowed and forbidden triad types and the sets of reduced (called selected) allowed and forbidden triad types are then used to generate networks with a given blockmodel structure. For a blockmodel type which cannot be generated successfully based only on the types of triads, some other local network structures are considered.

Beside the different types of triads, other subgraph types of a size higher than three can be used to generate networks with a given blockmodel. Milo et al. \cite{milo2002network} confirmed that different motifs are common in different empirical networks, where motifs are defined as "patterns of interconnections occurring in complex networks at numbers that are significantly higher than those in randomized networks". Different types of triads, rather than motifs, are considered chiefly because they are the smallest sociological unit from which the dynamic of a multi-person relationship can be observed \cite{davis1977clustering}. 
    
Various types of algorithms can be applied to generate networks with a given blockmodel, considering only different triad types. In this study, the Relocating Links algorithm (RL algorithm) and the Monte Carlo Multi Chain algorithm (MCMC algorithm) are used. If the generated structures with the selected set of triads, obtained using both algorithms, are very similar and close to the assumed ideal structure, one may conclude it is possible to generate networks with the assumed blockmodel structure by only considering the selected types of triads. On the other hand, if the generated networks are not similar and in line with the assumed blockmodel structure, one may speculate whether this is a consequence of the specifics of the algorithms or that the set of selected local structures is insufficient to generate this specific blockmodel.

In this study, it is assumed that the assignment of a unit to a group is unknown. Considering the information on the group assignment would require a different methodological approach. It is also assumed that the units' characteristics are not known. Kogut \cite{kogut2000network} reported that a certain structure's emergence in a network is often the consequence of rules that generate self-organisation dynamics. These rules do not need to be technological in origin, but can also reflect institutional or cultural norms and are also deeply embedded in the social identity of the units, meaning they are often invisible or unknown when examining an empirical social network. 
 
This paper is organised in the following way: first, the global network structures are described in terms of blockmodels and then the local network structures (namely different triad types) are presented. In this context, the sets of allowed and forbidden triad types and the sets of selected allowed and selected forbidden triad types are proposed for each blockmodel type. Considering these, random networks with the chosen blockmodels are generated using the RL algorithm and the MCMC algorithm. Generating a hierarchical blockmodel without complete blocks on the diagonal is further discussed in a separate section (Improvement of the hierarchical blockmodel without complete blocks on the diagonal) while the section Concluding remarks about generating networks with triads briefly summarises the ability to generate networks with a given blockmodel considering only different triad types. Some limitations and further research ideas are presented in the Conclusion. 

\section*{Global network structures}
    
Here, global network structures are defined by blockmodels. A blockmodel is a reduced network in which the units are clusters (positions) of units from the network. The term block refers to a submatrix showing the links between two clusters (positions) while the links in the network represent relationships between the positions \cite{doreian2005generalized}. Using the blockmodeling procedure, a blockmodel can be derived from a given empirical network. Blockmodeling is an approach for reducing a large, potentially incoherent network to a smaller, comprehensible and interpretable structure \cite{doreian2005generalized}. It can entail either a direct or indirect approach. 
    
 Several blockmodeling approaches have been developed to establish the best blockmodel structure according to the given network and equivalence \cite{doreian2005generalized}. In this study, pre-specified direct blockmodeling (generalised blockmodeling) was used. Here, the blockmodeling procedure is a local optimisation procedure. The solution is optimised with a relocation algorithm which minimises the value of the criterion function \cite{batagelj1997notes, batagelj1998fitting, doreian1994partitioning}. The criterion function reflects the difference between the ideal blockmodel and the empirical (current) solution. Compared to indirect blockmodeling, direct blockmodeling produces a solution with a lower criterion function value. Further, in the case of direct blockmodeling, the risk of a local optimum exists and therefore the algorithm must be repeated several times in the hope of obtaining the global optimum, and its computational complexity is high when a larger number of units is analysed.
    
There are several well-known and studied blockmodel types, e.g., cohesive, core-periphery, hierarchical and transitivity \cite{doreian2005generalized, wasserman1994social}. Even though all these structures have often been studied, there are different definitions of them \cite{borgatti2000models}. In this study, the definitions of blockmodel structures are taken from Doreian et al. \cite{doreian2005generalized}. 

\begin{figure}[!h]
\caption{{\bf Cohesive blockmodel with different level of errors.} 
(A) level of errors is 0 (ideal blockmodel), (B) level of errors is 0.25, (C) level of errors is 0.50, (D) level of errors is 1 (random network)}
\label{Fig2}
\includegraphics[height=1\textwidth, angle = 270]{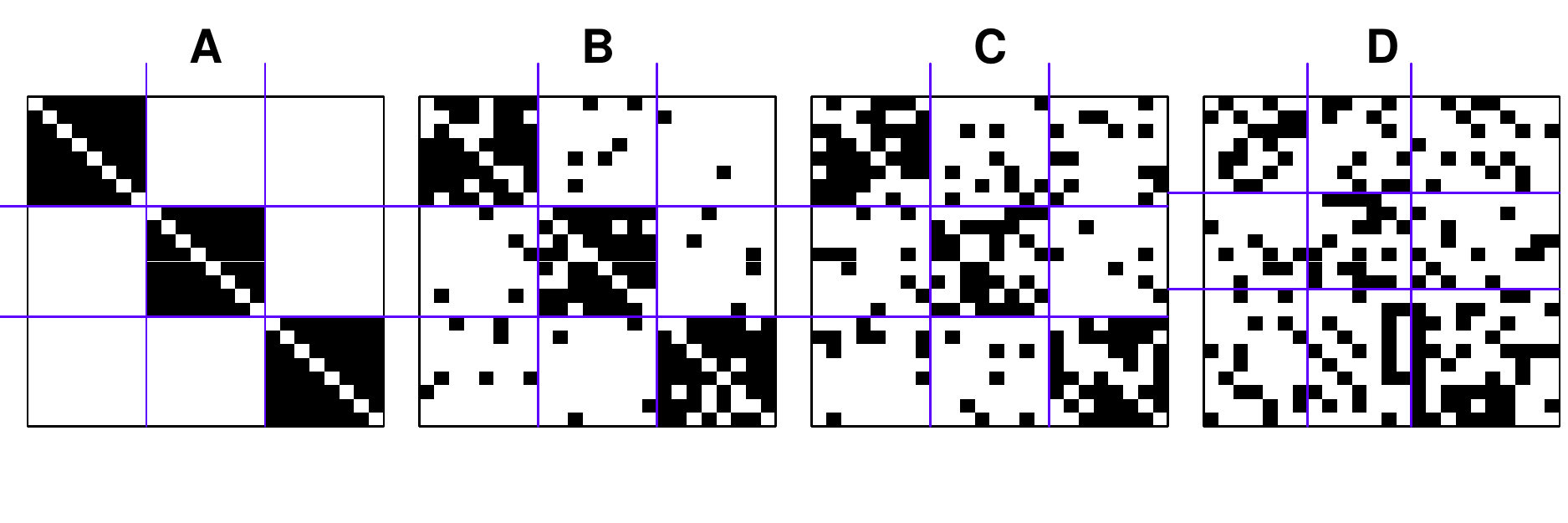}
\end{figure}

Structural equivalence \cite{lorrain1971structural} is considered in all blockmodels. It was shown \cite{doreian2005generalized} that in the case of structural equivalence only complete and null blocks exist. In ideal complete blocks, all possible links are present while in ideal null blocks no link exists. The left matrix in Fig~\ref{Fig2} represents a blockmodel with three groups (positions). There are only complete and null blocks where structural equivalence is completely satisfied . Such a blockmodel is a pre-specified blockmodel and called an "ideal blockmodel". In the blockmodel in Fig~\ref{Fig2}A, the diagonal blocks are complete and all the others are null blocks. This type of blockmodel is called a cohesive blockmodel. Usually, the empirical blockmodels being considered are not completely consistent with the selected equivalency and errors exist when comparing such a blockmodel with an ideal one. An error exists when a link is present in a null block or there is a non-link in a complete block. Fig~\ref{Fig2}B-D illustrates a cohesive blockmodel with different levels of errors. A level of errors is measured on a scale between 0 and 1, where 0 corresponds to an ideal blockmodel and 1 corresponds to a totally randomised network with the same density as in the ideal blockmodel. The level of errors increases linearly as ties are moved from complete blocks to null blocks, until the densities of both block types are the same (the level of errors then equals 1). In such a network, it is not possible to distinguish between blocks.
    
\begin{figure}[!h]
\caption{{\bf Different types of blockmodel structures: graph (left) and matrix (right) representation.} (A) cohesive, (B) symmetric core-periphery, (C) asymmetric core-periphery, (D) hierarchical without complete blocks on the diagonal, (E) hierarchical with complete blocks on the diagonal, (F) transitive without complete blocks on the diagonal, (G) transitive with complete blocks on the diagonal}
\label{Fig3}
\includegraphics[width=1\textwidth]{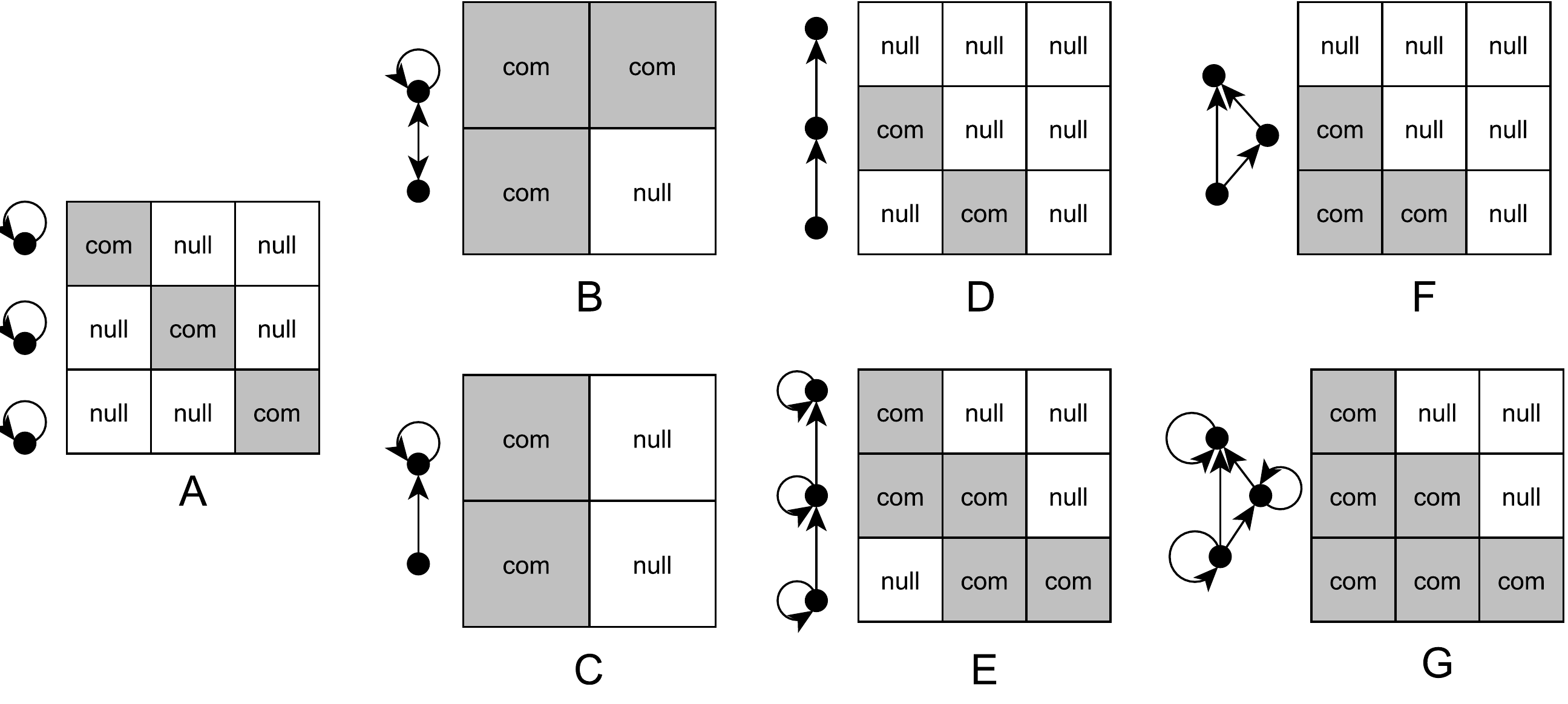}
\end{figure}

The following types of ideal blockmodels are defined and considered:
\begin{itemize}
\item	\textbf{Cohesive blockmodel} is visualised in Fig~\ref{Fig3}A. With this blockmodel, several internally highly connected clusters of units (positions) are present. The units from different clusters are not linked to each other. This is a very basic network structure type and was also studied, e.g. in the context of the structural organisation of the brain \cite{shen2015network}.
\item The most common \textbf{core-periphery blockmodel} consists of one group of units which are highly internally linked to each other. Peripheral units which are not linked to each other are also assumed in this type of network. The core-periphery blockmodel is called symmetric when the links between the peripheral and core units are mutual (Fig~\ref{Fig3}B) and asymmetric when only the peripheral units are linked to the core ones (Fig~\ref{Fig3}C) or when only the core units are linked to the peripheral ones. The core-periphery blockmodel structure lies in the middle of several extreme properties, e.g. clique vs. star configurations, network assortativity vs. network disassortativity, hierarchy vs. non-hierarchy, etc. \cite{csermely2013structure}. 

The core-periphery model is often associated with the existence of elites. An elite group is a small group of units that are all linked to each other (core). Compared to peripheral units, core units have greater prestige, usually defined by a higher number of incoming links (higher in-degree).

A clear core-periphery blockmodel was found among high school students where a link between students exists if the first student asked the second one to lend their study notes \cite{batagelj2004generalized}. It was also found when studying individual creative performances in the Hollywood film industry \cite{cattani2008core}, in the analysis of metabolic networks \cite{da2008centrality}, and in many studies of scientific co-authorships \cite{hu2008visual, cugmas2016stability, chinchilla2012blockmodeling, hu2008visual}.
\item	A \textbf{hierarchical blockmodel} consists of several groups of units which can be ordered into a hierarchy based on the direction of the links between the clusters. The units inside the groups can be either linked to each other (Fig~\ref{Fig3}D) or not (Fig~\ref{Fig3}E). A hierarchical structure is often associated with companies' organisational structure \cite{oberg2008hierarchical, ahuja1998network}.
\item A \textbf{transitivity blockmodel} (Fig~\ref{Fig3}F and Fig~\ref{Fig3}G) is similar to a hierarchical model. The only difference is that units from the groups on the lowest level are linked to all groups on the upper levels. This results in many transitivity relations (a relation $R$ on a set $A$ is called transitive if, for any $a,b,c \in A$ the conditions $aRb$ and $bRc$ imply $aRc$) which are very frequent when networks are formed among humans \cite{leinhardt1973development, schaefer2010fundamental} and animals \cite{chase1982dynamics, fararo1986state, skvoretz1996social}. In the literature, both hierarchical and transitive global network structures are often called hierarchical.   
\end{itemize}
    
\section*{Algorithms for generating networks}
    
To explain the impact of local mechanisms on global network structures or to characterise the global network structures in terms of local network structures, different statistical models have been developed \cite{toivonen2009comparative}. Many of these models capture different global network characteristics such as a specific distribution of in-degree or out-degree, the clustering coefficient or the small-world effect and less the specific global configuration of links in the network (e.g. in terms of the blockmodel structure) \cite{kejzar2007}.

To generate networks with a given blockmodel by considering different types of triads, two similar algorithms are used: the RL algorithm and the MCMC algorithm implemented in the \texttt{ergm} package implemented in R \cite{hunter2008ergm}. They both share the assumption that the units tend to create such a constellation of links that would result in a desirable distribution of subgraphs of size three or other characteristics in the network. Following the distinction between network evolution models (NEM), network attribute models (NAM) and ERGM \cite{toivonen2009comparative}, the RL algorithm can be classified in the NEM category. While NEM are primarily used to study how a specific rule (or set of rules) about creating and dissolving links affects the global network structure, the ERGM are used to check to what extent the global network structure can be explained when considering the structure of links and/or characteristics of the units. It can also be used to generate networks based on estimated or fixed parameter values. Both approaches are described and compared in more detail in the following sections.

\subsection*{Generating networks with the Relocating Links algorithm (RL algorithm)}
    
The RL algorithm (see Algorithm~\ref{RLalgorithm}) assumes that all considered local network statistics for the case of an ideal network are represented by the vector $\mathfrak{T}$. The number of elements $g$ of this vector equals the number of local network statistics considered. The numbers of different types of triads are considered here, but some other local network statistics could also be chosen. The distribution of all or only a subset of all triad types can be given (for forbidden triad types, corresponding values of $\mathfrak{T}$ equal zero). Beside $\mathfrak{T}$, the initial random network $Y_r$ has to be given. Before the iterative procedure starts, the $Y_r$ is saved as a new network $Y_{new}$.
    
The iterative procedure is repeated many times. Upon each iteration, a pair of linked units $i$ and $j$ and a pair of unlinked units $k$ and $l$ are randomly chosen. Then, the link between $i$ and $j$ is dissolved and the link between $k$ and $l$ is established. The modified network is saved as $Y_p$ (the proposed network). From the proposed network $Y_p$, the number of each triad type considered $\mathfrak{T}_p$ is calculated. The proposed network is saved as $Y_{new}$ (the new network) if the CR ratio is greater than 1. The CR ratio is defined as 
    
\begin{eqnarray}
\label{eq:cr}
CR = \frac{\sum_{i=1}^g {\Big (}(\mathfrak{T}_p - \mathfrak{T})^2{\Big )}_i} 
{\sum_{i=1}^g {\Big (}(\mathfrak{T}_{new} - \mathfrak{T})^2{\Big )}_i}
\end{eqnarray}
    
Then, the new iteration is performed and, after many iterations, the last $Y_{new}$ is the final solution. Besides the $Y_{new}$, the values of $CR$ can be saved and further analysed. 
    
\begin{algorithm}
\caption{The Relocating Links algorithm}\label{RLalgorithm}
\begin{algorithmic}[1]
    \Require $\mathfrak{T}$ \Comment{$\mathfrak{T}$ denotes the distribution of local network statistics in an ideal network}
    \Require $Y_r$ \Comment{$Y_r$ denotes a random network}
    \Require $k$ \Comment{$k$ denotes the number of iterations}
    \State $Y_{new} \gets Y_r$
      \State $Y_{p} \gets Y_r$
      \For{$k$ in $1:k$}
    \State randomly select a tie $y_{i,j}$ in $Y_{new}$
      \State randomly select a non-tie $y_{k,l}$ in $Y_{new}$
      \State transform a tie $y_{i,j}$ to a non-tie in $Y_p$
      \State transform a non-tie $y_{k,l}$ to a tie in $Y_p$
      \If{$CR > 1$}  \Comment{$CR$ is defined in Eq~\ref{eq:cr}}
    \State $Y_{new} \gets Y_p$ 
      \Else \State $Y_p \gets Y_{new}$
      \EndIf
    \EndFor \\
    \Return $Y_{new}$
\end{algorithmic}
\end{algorithm}
    
Compared to the MCMC algorithm introduced in the next section, the RL algorithm is more deterministic since a link is only allocated if the distribution of the triads of the proposed network is closer to the distribution of the triads in the case of an ideal blockmodel. This may result in lower variability of the global network structure of generated networks when the RL is used since, compared to the MCMC algorithm, RL strive to generate networks with the exact number of the selected types of triads. However, the risk of a local optimum exists which could be avoided by further improving the algorithm. Moreover, RL is computationally very intensive: as will be illustrated later, a higher number of iterations is required, especially in the case of denser networks. 
    
\subsection*{Generating networks with the MCMC algorithm}
    
To describe how the networks were generated using the MCMC algorithm, Exponential Random Graph Modelling (ERGM) has to be defined. Let us consider a random network $Y$ ($y$ is a given empirical network) consisting of $N$ units. Here, the link between the $i$-th and the $j$-th unit can be represented by a random variable $Y_{ij}$, while the set of all possible random networks of this size is denoted by $\mathcal{Y}$. The distribution of $Y$ can be written as
    
\begin{eqnarray}
\label{eq:ergm}
P_{\theta, \mathcal{Y}} (Y=y) = \frac{exp\{\theta^T g(y)\}}{\kappa (\theta, \mathcal{Y})}
\end{eqnarray}
    
\noindent where $y \in \mathcal{Y}$. Here, $\theta$ is a vector of coefficients while $g(y)$ is the vector of statistics obtained for $y$. The normalisation constant $\kappa (\theta, \mathcal{Y})$ in the numerator is needed to ensure the sum of probabilities equals 1. 
    
Different methods can be used \cite{corander1998maximum, strauss1990pseudolikelihood, hyvarinen2007connections} to estimate the parameters $\theta$. After that, one can generate random networks based on the model obtained. To do this, several types of MCMC algorithms have been proposed. Generally, the start is represented by a network in $\mathcal{Y}$. Then, based on a uniform distribution one of the links or non-links is chosen. According to the model, the probability of establishing or dissolving a link is calculated and then, based on this probability, the chosen link or non-link is established or dissolved. The process is iterative. For each iteration, the change statistic is calculated, namely the change in the values of the estimated statistics before and after the change in the link between $i$ and $j$. The iterative process stops when approximate convergence to $P_{\theta_0, Y} (Y=y)$ is reached \cite{hunter2008ergm}. There are several versions of the described algorithm which chiefly differ regarding how the probability of establishing or dissolving a link is calculated. In the used \texttt{ergm} package, the Metropolis-Hastings algorithm is implemented and used. 
    
A very common problem with the MCMC algorithm is the multimodal distribution of sufficient statistics \cite{snijders2002markov, jonasson1999random}. A so-called degenerative model emerges when the model poorly fits the empirical data (e.g. due to inappropriately chosen terms), resulting in generated networks that do not fit the empirical network. These generated networks are often without any link or all possible links \cite{handcock2003assessing}. One possible solution is to restrict the class of networks considered to be possible under the model by fixing the number of links. 
    
As described, the method most often used to estimate the parameters is MCMC-MLE which can sometimes be computationally hard to estimate. In our study, the parameters can be estimated based on networks with a given blockmodel without or with only very low levels of errors. Using this approach, the estimation algorithm does not converge in many cases, probably due to the high level of multicollinearity of the triads. In addition, from the end-user point of view, estimating the values of all parameters for each blockmodel type would be very hard.
    
Instead, the values of the ERGM parameters $\theta$ are arbitrarily set to 2 (allowed) or -2 (forbidden). It has been clearly shown that some types of triads are much more likely to appear in an ideal network (compared to a random network). By setting all the parameters' values to 2 or -2, we essentially assume that all types of allowed triads have the same importance (and similar for all forbidden types of triads). Such a setting is critical when all types of triads are included in the model and result in a relatively unstable model, particularly when the density is not fixed.

All types of triads are generated considering the number of links fixed (to the same value as in ideal networks) on one hand, and free (with the density being the variable) on the other. In the case of the latter, the value of parameter \texttt{edge} is set to such a value that the mean density of 30 generated networks lies within the interval of ideal density $\pm 0.05$.

\section*{Choosing triads for different types of ideal blockmodels}
  
When generating networks with a specific type of a blockmodel, according to different triad types, all triad types or only a subset of all triad types can be considered. This is particularly important when generating networks with the RL algorithm where the distribution of triads has to be known in advance for each type of ideal blockmodel separately. Here, it has to be pointed out that the distributions of triads can vary among the same type of blockmodel with a different number of positions.

Since the number of different triads is affected by the network density \cite{faust2006comparing}, the value of the A-measure (the ratio between the absolute number of a certain type of triad in an ideal blockmodel and the mean number of such triads in the totally randomised network of the same density -- see Appendix S23 for more information about generating totally randomised networks and the networks with a given level of errors) can be used to minimise the number of different triad types needed to generate the networks with a selected blockmodel. 

The classifications of allowed and forbidden triad types for different blockmodel types are presented in the next section followed by the classifications of selected allowed and selected forbidden triad types, based on values of the A-measure.

\subsection*{Allowed and forbidden triad types}
    
The triad types can be classified in the set of allowed or in the set of forbidden triad types for each blockmodel type, based on the counts of triad types in an ideal blockmodel. Triad types with the count equal to zero are said to be forbidden in a given blockmodel and are thus classified in the set of forbidden triad types (for a given blockmodel). All the other triad types are classified in the set of allowed triad types. This classification is essential for the MCMC algorithm as it determines the values of the appropriate parameters in the ERGM model (see the previous section).
    
A more detailed insight into how common a certain triad type is for a certain blockmodel can be obtained by interpreting the A-measure values. The A-measure allows the relative number of triads to be compared within a certain type of blockmodel and also the relative number of triads between different types of blockmodels. In order to obtain values of the A-measure, 10,000 totally random networks for each ideal blockmodel type were generated. The A-measure values are presented in Table~\ref{selected_terms}. Values greater than 1 indicate triad types that are more likely to occur in an ideal network structure than would be expected in randomised networks. Such are complete subgraphs of size three (a triad of type 300) in a cohesive blockmodel. 
    
When the A-measure value is close to 1, the number of triads in the case of an ideal network structure is close to the number of triads in totally randomised networks. This could be an indicator that their occurrence is mainly a consequence of the density rather than the type of blockmodel. The A-measure values in the empty cells in Table~\ref{selected_terms} equal zero and therefore denote forbidden triads. It can be seen that the values of the A-measure of certain triad types exceed zero in some but not all blockmodel types. 
    
Reducing the number of triad types used to generate networks with a given blockmodel can be beneficial in several ways. For example, it can help to identify the main (e.g. social) mechanisms that cause a given blockmodel structure to be formed.

In addition, there are practical reasons which differ with respect to the algorithm used. For the RL algorithm, the reduction to only forbidden triad types (or a subset of forbidden triad types) is especially appealing as it does not require knowledge of the exact distribution of triad types in the ideal network (as this algorithm otherwise requires) because the frequency of all forbidden triad types is 0. 

The frequencies of different forbidden triad types are also not affected by the sizes and number of clusters. This means that, when generating networks by considering only the forbidden triad types, this information is not taken into account, which may be either desired or not. On the other hand, the frequencies of all allowed triad types contain all the information that is included in all (allowed and forbidden) triad types. Therefore, considering only all allowed triad types is hypothesised to be equivalent to considering all possible triad types. 

\subsection*{Selecting subsets of triad types}
    
Networks generated by the RL algorithm can still differ as the CR (see Eq.~\ref{eq:cr}) is computed slightly differently. For the MCMC algorithm, these issues are not relevant since the exact distribution of triad types is never taken into account when setting the parameter values. However, the MCMC algorithm is affected by multicollinearity, which can be reduced by selecting only certain triad types. Given the point of this algorithm, it is best to select only a small number of relatively different triad types. 

The sets of allowed and forbidden triad types can be further reduced to the selected allowed and selected forbidden triad types. This means that not all possible or all allowed or all forbidden types of triads are considered. One could choose the most appropriate triad types based on observations of the A-measure in the case of networks without errors as was done in the previous subsection. Yet, to obtain a better selection of triad types, it is beneficial to observe the A-values for networks with different levels of errors (10,000 networks with the value of the level of errors from 0.2 to 1 with step 0.2 are generated) as done in this study (see Appendix S22 for a more detailed description).

The selected triad types are shown in grey in Table~\ref{Tab1}. It may be seen that only a few triad types are allowed for each type of blockmodel. Almost all of these types of triads are chosen for all types of blockmodels. The exceptions are triad type 021C in the case of a hierarchical blockmodel without complete blocks on the diagonal and triad type 300 in the case of a transitivity blockmodel with complete blocks on the diagonal. On the other hand, for some blockmodel types only a small number of all forbidden triad types is selected, e.g. in the case of an asymmetric core-periphery only one, and in the case of a cohesive blockmodel only  two. 
    
\begin{table}[!h]
\begin{adjustwidth}{-2.25in}{0in}
\caption{{\bf A-measure values and the classification of allowed and forbidden triad types, and selected allowed and forbidden triad types for different types of blockmodels.} Values greater than zero denote allowed types of triads while the values which equals zero (empty cells) denotes forbidden types of triads; grey color denotes selected triad types.}
\label{Tab1}
\footnotesize 
\begin{tabular}{|p{1cm}|p{1.8cm}|p{2cm}|p{2cm}|p{2.2cm}|p{2.2cm}|p{2cm}|p{2cm}|}
\hline
    & COHESIVE & ASYMMETRIC  CORE-PERIPHERY&SYMMETRIC  CORE-PERIPHERY&HIERARCHICAL WITHOUT COMPLETE BLOCKS ON THE DIAGONAL&HIERARCHICAL WITH COMPLETE BLOCKS ON THE DIAGONAL &TRANSITVITIY WITHOUT COMPLETE BLOCKS ON THE DIAGONAL& TRANSITIVITY WITH COMPLETE BLOCKS ON THE DIAGONAL\\ \hline
    003 &2.3&7.1&7.2&1.5&&1.1&~\\ \hline
    300  & 96.3\cellcolor{gray!25}&7.5\cellcolor{gray!25}&2.7\cellcolor{gray!25}&&3.7\cellcolor{gray!25}&&1.2\\ \hline
    120D &&8.2\cellcolor{gray!25}&\cellcolor{gray!25}&\cellcolor{gray!25}&4.1\cellcolor{gray!25}&&5.1\cellcolor{gray!25}\\ \hline
    120U &&&\cellcolor{gray!25}&&4.1\cellcolor{gray!25}&&5.1\cellcolor{gray!25}\\ \hline
    102  &10.2\cellcolor{gray!25}&\cellcolor{gray!25}&\cellcolor{gray!25}&&5.8\cellcolor{gray!25}&&\cellcolor{gray!25}\\ \hline
    021C &&&\cellcolor{gray!25}&2.2&3.1\cellcolor{gray!25}&\cellcolor{gray!25}&\cellcolor{gray!25}\\ \hline
    021U &\cellcolor{gray!25}&8.2\cellcolor{gray!25}&\cellcolor{gray!25}&4.0\cellcolor{gray!25}&\cellcolor{gray!25}&5.1\cellcolor{gray!25}&\\ \hline
    021D &\cellcolor{gray!25}&&\cellcolor{gray!25}&4.0\cellcolor{gray!25}&\cellcolor{gray!25}&5.1\cellcolor{gray!25}&\\ \hline
    030T &&&&&&3.5\cellcolor{gray!25}&3.5\cellcolor{gray!25}\\ \hline
    201  &&&6.6\cellcolor{gray!25}&\cellcolor{gray!25}&\cellcolor{gray!25}&&\cellcolor{gray!25}\\ \hline
    120C &&&\cellcolor{gray!25}&&\cellcolor{gray!25}&\cellcolor{gray!25}&\cellcolor{gray!25}\\ \hline
    111D &&&&\cellcolor{gray!25}&&\cellcolor{gray!25}&\cellcolor{gray!25}\\ \hline
    111U &&&&&&\cellcolor{gray!25}&\cellcolor{gray!25}\\ \hline
    030C &&&&&&&\\ \hline
    210  &&&&&&&\\ \hline
    012  &&&&&&&\\ \hline
\end{tabular}
\end{adjustwidth}
\end{table}
    
\section*{Simulation design}
    
The whole simulation process is accomplished by using the R programming language \cite{team2000r}. Using each method (RL algorithm and MCMC algorithm), $k=50$ networks (of size $n=24$ units) with a given blockmodel structure are generated for each selected set of triads. Each generated network is randomised. Pre-specified blockmodeling (also called generalised blockmodeling) is applied to model networks and randomised networks where the number of clusters is set as in the ideal networks (to two or three clusters) (the \texttt{blockmodeling} \cite{ziberna2008blockmodeling} package implemented in R). For each generated network, the value of the blockmodeling criterion function is calculated. 

Here, it should be highlighted that there may be bias in the values of the criterion function, where the networks are generated by the RL algorithm and all allowed triad types are considered. This is because the information on the number and sizes of the groups is embedded in the frequencies of the different allowed triad types when using the RL algorithm. Yet this is not the case when the MCMC algorithm is used and/or other subsets of triads are considered.

However, as the criterion function is not generally comparable for different blockmodels, the Mean Improvement Value (MIV) is calculated as

\begin{eqnarray}
\label{eq:miv}
MIV = 1 - \frac{1}{k} \sum_{i=1}^k {\Big (}\frac{P^m_i}{P^r_i} {\Big )}
\end{eqnarray}
    
\noindent where $P^r_i$ is the value of the criterion function of the $i$-th randomised network and $P^m_i$ is the value of the criterion function of $i$-th network generated based on the model. The MIV obtained on a network with a different number of units is generally not comparable.

For each type of generated network (the RL algorithm, the MCMC algorithm with fixed density, and the MCMC algorithm with non-fixed density), visualisations of the $P^r$ and the $P^m$ are presented in Appendix~S19-S21 for each type of blockmodel and each combination of triad types considered. The corresponding values of the MIV are visualised in Fig~\ref{Fig5}, Fig~\ref{Fig7} and Fig~\ref{Fig8}.

\section*{Results}
    
This section is organised in several subsections. First, the global network structures of the networks generated with the RL algorithm and MCMC algorithm (fixed and non-fixed density) are evaluated. For each algorithm, the networks generated by considering different sets of triad types are compared. Then, considerations are presented of certain other local network structures in the evant triads do not generate the networks with the expected blockmodel. Finally, a general statement on generating networks using triads is given.
    
\subsection*{Networks generated with the RL algorithm}
    
When all triad types are considered, the overall MIV is around 72~\%, which is more than for all other sets of triads considered (Fig~\ref{Fig5}). On the other hand, the MIV corresponding to networks generated based only on all forbidden triads or all allowed triads is slightly lower or the same. What is outstanding is the symmetric core-periphery with the lowest MIV varying between 11 and 32~\% among the different models (all, all allowed, or all forbidden triad types). As has been emphasised, when the network is very dense the RL algorithm is less effective in finding the right link to relocate. This is expressed in a very small peripheral part in the case of a symmetric core-periphery blockmodel.
    
The MIVs are usually lower when all forbidden triad types are considered. The MIVs corresponding to the cohesive blockmodel are very similar, yet the structure of the blockmodels so generated is different when only the set of forbidden triad types is considered (the cluster sizes are more variable).

When comparing the different blockmodel types, the highest MIV is observed in the case of an asymmetric core-periphery blockmodel (98~\% when all allowed or all forbidden types of triads are considered), in the case of a transitivity blockmodel without complete blocks on the diagonal (95~\% when all allowed types of triads are considered and 94~\% when all types of triads are considered) and in the case of a transitivity blockmodel with complete blocks on the diagonal (94~\% when all allowed triad types are considered and 92~\% when all triad types are considered). In the latter case, there is quite considerable variability in the cluster sizes when all forbidden types of triads are considered. To be more precise, the tendency to form one cluster with a relatively high number of units and two clusters with a lower number of units is present. This happens because different types of triads can be present in different parts of the network. 

\begin{figure}[H]
\caption{{\bf Some generated networks with the expected hierarchical blockmodel without complete blocks on the diagonal.} The RL algorithm is used and all triad types are considered.}
\label{Fig4}
\includegraphics[height=1\textwidth, angle = 270]{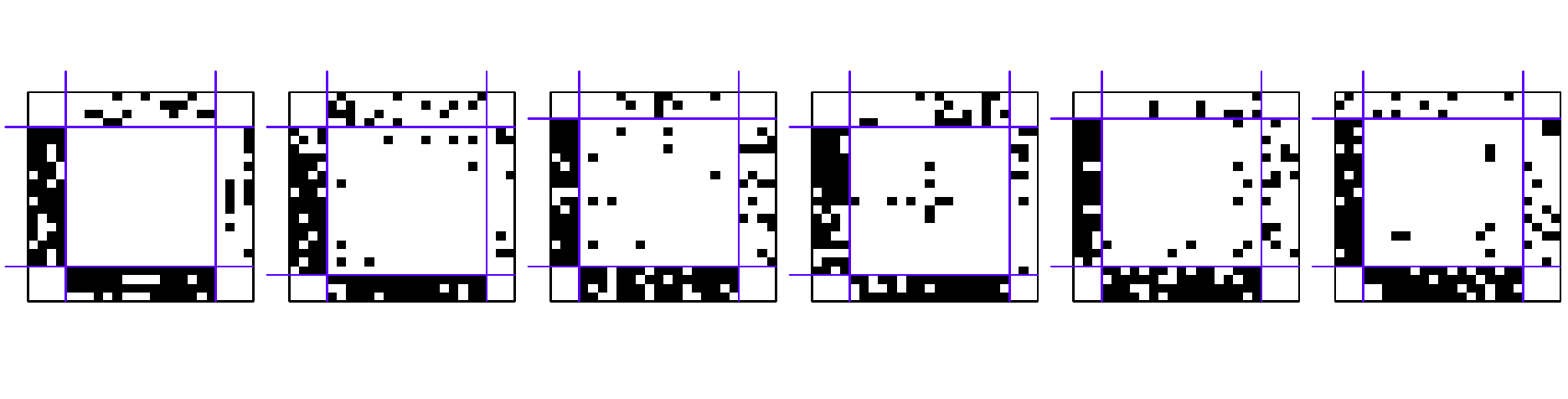}
\end{figure}

\begin{figure}[H]
\begin{adjustwidth}{-2.25in}{0in}
\caption{{\bf The mean MIV for each blockmodel type (generated by the RL algorithm) and selected set of triad types.} (A) cohesive; (B) symmetric core-periphery; (C) asymmetric core-periphery; (D) hierarchical without complete blocks on the diagonal (E) hierarchical with complete blocks on the diagonal; (F) transitivity without complete blocks on the diagonal; (G) transitivity with complete blocks on the diagonal, (1) allowed and forbidden triad types, (2) allowed triad types, (3) forbidden triad types. Note: only networks of transitivity with complete blocks on the diagonal blockmodel type and hierarchical with complete blocks on the diagonal blockmodel type were generated by considering the selected allowed triad types.}
\label{Fig5}
\includegraphics{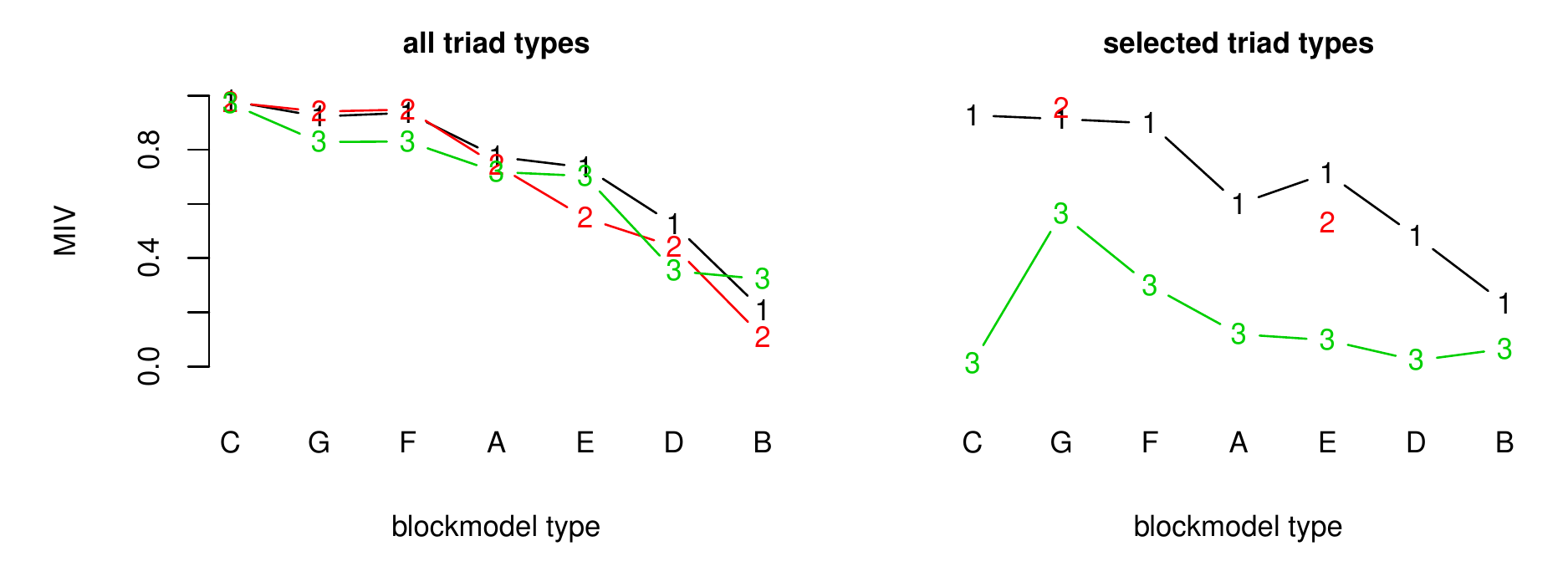}
\end{adjustwidth}
\end{figure}

When generating networks with a hierarchical blockmodel without complete blocks on the diagonal, a blockmodel structure which is not assumed emerges. Instead, there are links in the blocks below the diagonal of the matrix and in the blocks above the diagonal. This means there are links from the top to the lowest clusters and from the low clusters to the higher clusters. On the level of units, only asymmetric links are possible. However, the density is still higher in complete than in null blocks (Fig~\ref{Fig4}), which may be a consequence  of the optimisation algorithm for pre-specified blockmodeling.  
    
When all allowed triads are included in the process of generating networks, one would expect a similar MIV as when all triads are included in the model. This is because all the information for generating the networks embedded in all triad types is also embedded in only allowed triad types (as all the rest have a count of 0).

The set of all allowed triad types and the set of triads with selected allowed types of triads vary only in the case of a hierarchical blockmodel with complete blocks on the diagonal and a transitivity blockmodel with complete blocks on the diagonal. The selection of triad types slightly improves the MIV in the case of both blockmodel types. In the former case, the blockmodel structure can be visually recognised in most, but not all, generated networks. On the other hand, there are very low levels of errors in all generated networks with a transitive blockmodel with complete blocks on the diagonal.

Comparing the networks generated with all forbidden triad types and the networks generated with only the selected forbidden triad type, the MIV is generally lower in the latter case for all types of blockmodels. By visually observing some generated networks, it is hard to recognise the assumed blockmodel structure, except for some transitive blockmodels with complete blocks on the diagonal.

\subsection*{Networks generated with the MCMC algorithm: fixed density}
    
Since the RL algorithm is more deterministic, it generally performs better than the MCMC algorithm. But when networks are denser the MCMC algorithm might perform better as this is the case when e.g., considering the set of all allowed types of triads when generating a symmetric core-periphery blockmodel. This is another reason one has to consider different algorithms when studying micro structures in the context of various global network structures using simulations. 

\begin{figure}[H]
\caption{{\bf Some empirical generated networks using the RL algorithm by considering all triad types.} By rows: (A) cohesive; (B) symmetric core-periphery; (C) asymmetric core-periphery; (D) hierarchical without complete blocks on the diagonal; (E) hierarchical with complete blocks on the diagonal; (F) transitivity without complete blocks on the diagonal; (G) transitivity with complete blocks on the diagonal.}
\label{Fig6}
\includegraphics[width=1\textwidth]{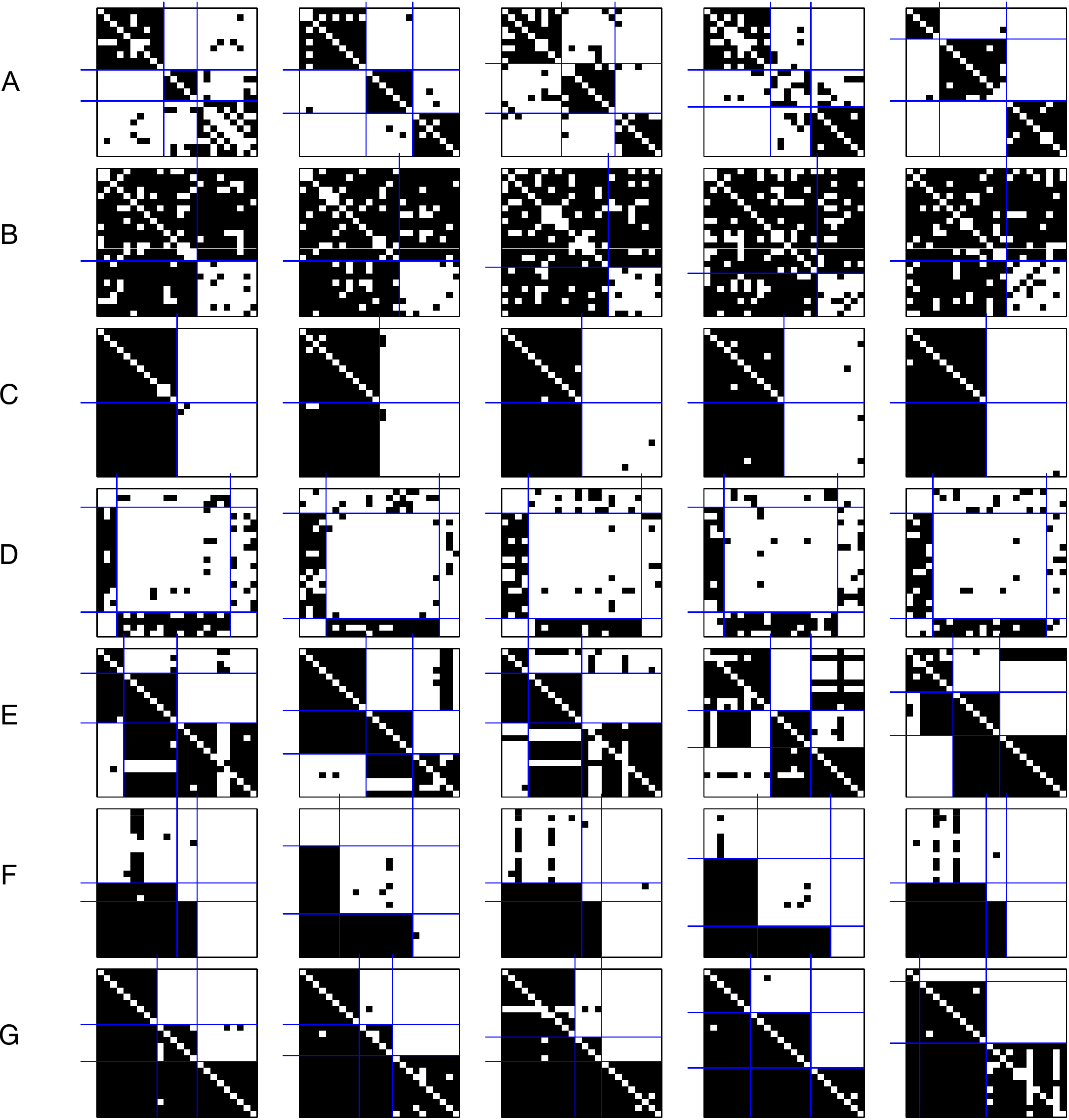}
\end{figure}

\begin{figure}[H]
\begin{adjustwidth}{-2.25in}{0in}
\caption{{\bf The mean MIV for each blockmodel type (generated by the MCMC algorithm with fixed density) and selected set of triad types.} (A) cohesive; (B) symmetric core-periphery; (C) asymmetric core-periphery; (D) hierarchical without complete blocks on the diagonal (E) hierarchical with complete blocks on the diagonal; (F) transitivity without complete blocks on the diagonal; (G) transitivity with complete blocks on the diagonal, (1) allowed and forbidden triad types, (2) allowed triad types, (3) forbidden triad types. Note: only the networks of transitivity with complete blocks on the diagonal blockmodel type and hierarchical with complete blocks on the diagonal blockmodel type were generated by considering the selected allowed triad types.}
\label{Fig7}
\includegraphics{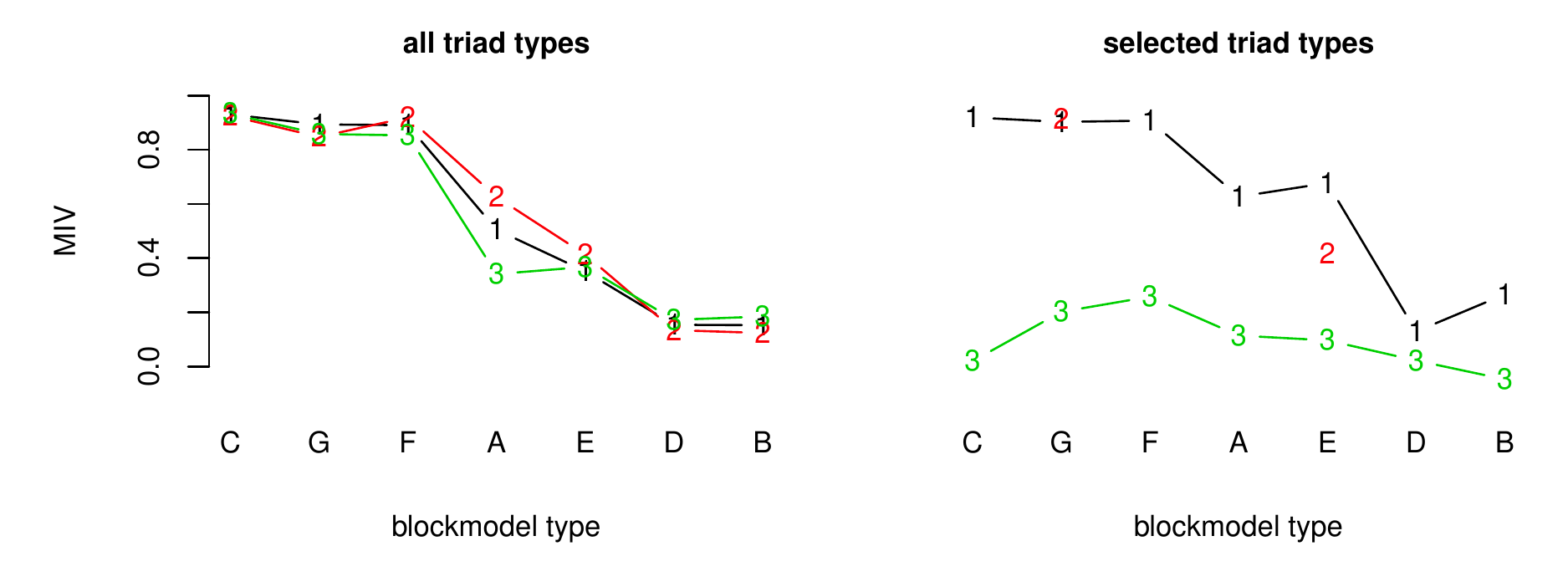}
\end{adjustwidth}
\end{figure}

When all possible triad types are considered, the overall mean MIV among all blockmodel types is higher when the networks are generated using the RL algorithm and lower when the networks are generated using the MCMC algorithm with a fixed density (Fig~\ref{Fig7}). Yet generated networks have an assumed blockmodel structure (Fig~\ref{Fig6}) with a relatively low level of errors, except the hierarchical one without complete blocks on the diagonal where the global network structure obtained is similar to that produced with the RL algorithm (considering selected allowed triad types) (see Fig~\ref{Fig4}). Further, the hierarchical structures with complete blocks on the diagonal and the cohesive one are less clear than the others.

Considering only all allowed or only all forbidden triad types does not produce networks with a significantly higher level of errors. The MIVs are lower when selected forbidden triad types are considered compared to the case when all forbidden triad types are considered. In this instance, the generated networks do not have the expected blockmodel. 

Further selection of the different types of triads that are allowed does not improve a hierarchical blockmodel with complete blocks on the diagonal, even though some MIVs indicate the opposite. Conversely, a further selection of allowed triad types improves the global structure of networks with an expected transitivity blockmodel with complete blocks on the diagonal.  
 
The further selection of all possible triad types (allowed and forbidden) improves all the MIV values, especially those corresponding to the hierarchical blockmodel with complete blocks on the diagonal and the cohesive blockmodel.

The generated networks with the expected hierarchical blockmodel structure without complete blocks on the diagonal are not in line with the expected global network structure. This is true for any set of triad types considered.

\subsection*{Networks generated with the MCMC algorithm: non-fixed density}
    
\begin{figure}[H]
\begin{adjustwidth}{-2.25in}{0in}
\caption{{\bf The mean MIV for each blockmodel type (generated by the MCMC algorithm with variable density) and selected set of triad types.} (A) cohesive; (B) symmetric core-periphery; (C) asymmetric core-periphery; (D) hierarchical without complete blocks on the diagonal (E) hierarchical with complete blocks on the diagonal; (F) transitivity without complete blocks on the diagonal; (G) transitivity with complete blocks on the diagonal, (1) allowed and forbidden triad types, (2) allowed triad types, (3) forbidden triad types. Note: only the networks of transitivity with complete blocks on the diagonal blockmodel type and hierarchical with complete blocks on the diagonal blockmodel type were generated by considering the selected allowed triad types.}
\label{Fig8}
\includegraphics{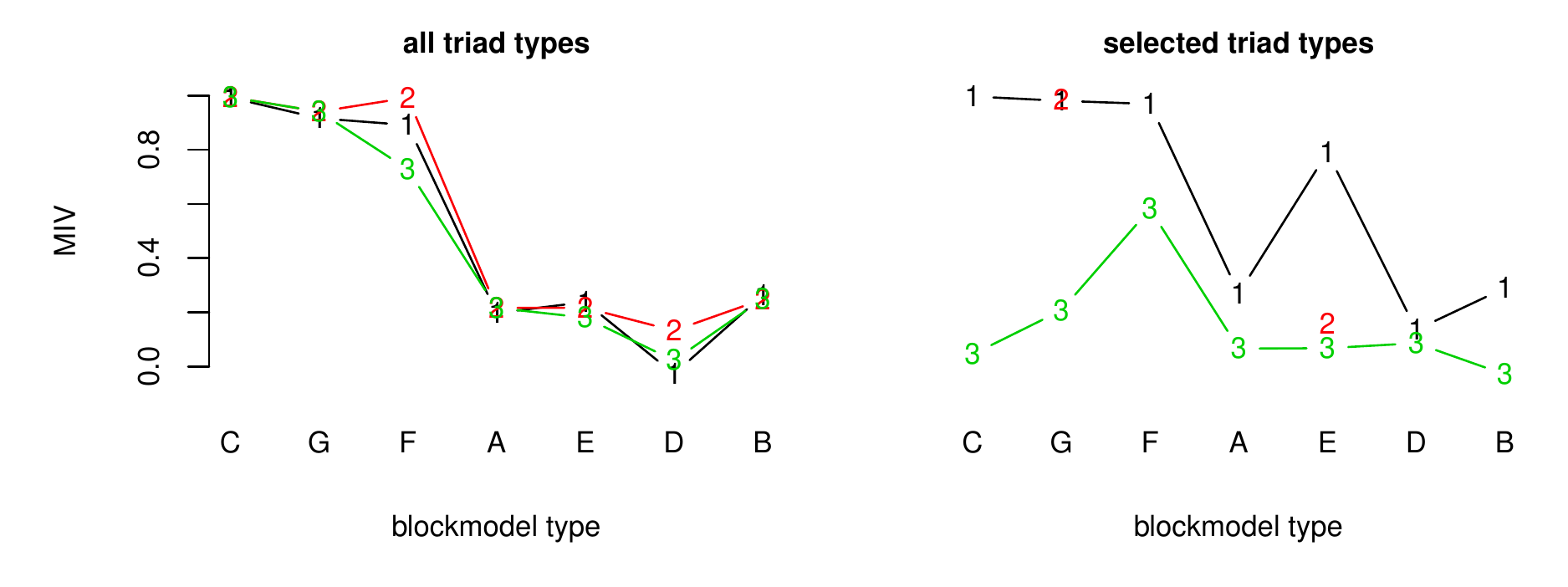}
\end{adjustwidth}
\end{figure}

In the event the initial networks are totally randomised ideal networks, networks generated using the MCMC algorithm with a non-fixed density are close to the networks with a fixed density (Fig~\ref{Fig8}). The further selection of different triad types is not seen as so important in the case of an asymmetric core-periphery blockmodel and a transitivity blockmodel (with or without complete blocks on the diagonal), while it improves the structure of generated networks with an expected symmetric core-periphery blockmodel and a hierarchical blockmodel with complete blocks on the diagonal. 
    
Here, it is noted that the way the initial networks are chosen has a great impact on the networks generated. In the case of the MCMC algorithm with a non-fixed density, considering the random networks (as initial networks) with the expected (the actual number becomes a random variable) number of links being equal to the number of nodes usually produces a very high number of totally empty or full generated networks. This is especially when all triad types are included in the model. In this study, the randomised ideal networks are used as initial networks, meaning the density of the initial networks is not variable and is the same as in the ideal networks.
    
\subsection*{Improvement of the hierarchical blockmodel without complete blocks on the diagonal}
    
\begin{figure}[H]
\caption{{\bf Some generated networks with a hierarchical blockmodel structure without complete blocks on the diagonal using different models.} The networks are generated by the MCMC algorithm with a non-fixed density. (a) All selected types of triads; (b) All selected types of triads and paths of length three; (c) All selected types of triads (with a higher parameter value for the triad of type 021D) and paths of length three}
\label{Fig9}
\includegraphics[height=1\textwidth, angle = 270]{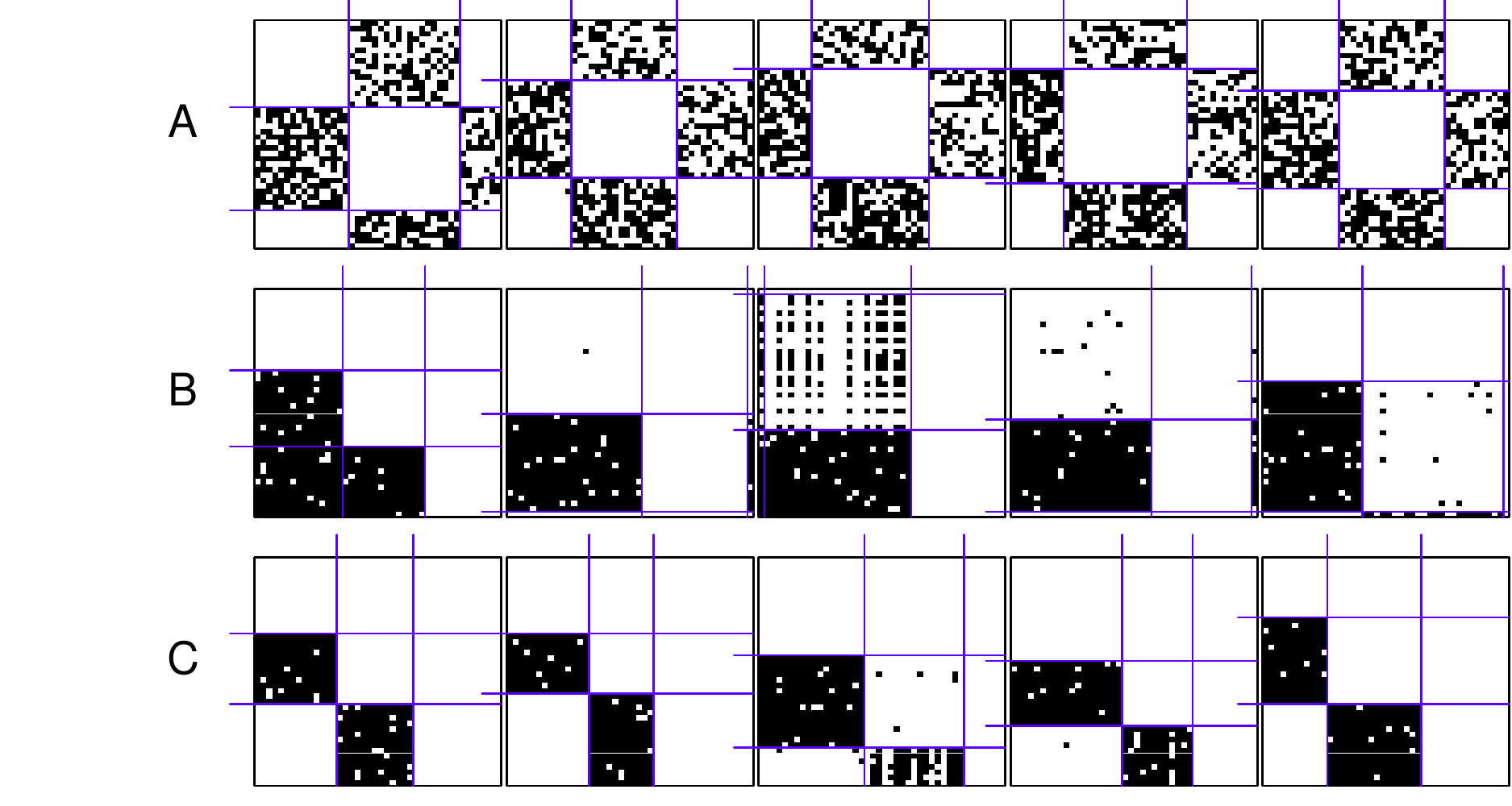}
\end{figure}

The proposed models for generating networks with a hierarchical blockmodel structure without complete blocks on the diagonal perform poorly. This is seen by the mean improvement values and the distribution of the values of the criterion function (see Appendix~S1-S21, Fig~\ref{Fig5}, Fig~\ref{Fig7}, Fig~\ref{Fig8} and Fig~\ref{Fig9}A for some empirical examples). 

The obtained blockmodel structure is often hierarchical but has additional links from the upper to the lower clusters and with all asymmetric links. This is especially typical of networks generated using the MCMC algorithm. Therefore, the main focus is put on the networks generated using the MCMC algorithm with a non-fixed density. The resulting global structure probably emerges since all considered triad types appear in all parts of the network. Their combination produces a network that is highly determined by paths of length three (e.g., $1 \rightarrow 2  \rightarrow 3  \rightarrow 2$, where digits denote clusters). The MIV is 0.17.

Therefore, omitting the links from the upper to the lower positions means the paths of length three are considered. Here, it should be pointed out that the number of triads is unit-based while the number of paths of length three is an edge-based count. However, an additional parameter paths of length three (in the case of networks with a different number of positions, paths of different lengths should be considered) is added to the model with the value of -2 (as forbidden). Networks generated using this model have the expected hierarchical structure but with only two groups (Fig~\ref{Fig9}B). From time to time, networks with a transitivity blockmodel without complete blocks on the diagonal are also produced (MIV=0.74).

To obtain three positions (instead of two), the parameter's value of triad type 021C has to be increased, e.g. to the value to 4. Such a model produces networks with a very clear hierarchical structure without complete blocks on the diagonal (Fig~\ref{Fig9}C). There are no errors in all generated networks in null blocks while some appear in complete blocks (MIV=0.93). 

All of the described networks are generated using the MCMC algorithm. When the RL algorithm is used, all allowed types of triads and paths of length three can be considered. In that case, some errors appear in both null and complete blocks, which is a consequence of the fixed density. However, with a higher number of iterations, the number of errors could also be lower. 

\subsection*{Concluding remarks about generating networks with triads}
    
Three approaches are proposed for generating networks with a given blockmodel structure. The networks so generated are compared with totally random networks of the same density. As expected, there are fewer errors in these networks when generated using the RL algorithm compared to the MCMC algorithm. With both approaches, the selection of triad types does not necessarily result in generated networks with higher or lower levels of errors.

\begin{figure}[H]
\begin{adjustwidth}{-2.25in}{0in}
\caption{{\bf Some examples of generated networks for each type of a blockmodel.} The networks are generated using the RL algorithm by selected allowed types of triads. (A) cohesive; (B) symmetric core-periphery; (C) asymmetric core-periphery; (D) hierarchical without complete blocks on the diagonal (E) hierarchical with complete blocks on the diagonal; (F) transitivity without complete blocks on the diagonal; (G) transitivity with complete blocks on the diagonal.}
\label{Fig10}
\includegraphics[angle = 270]{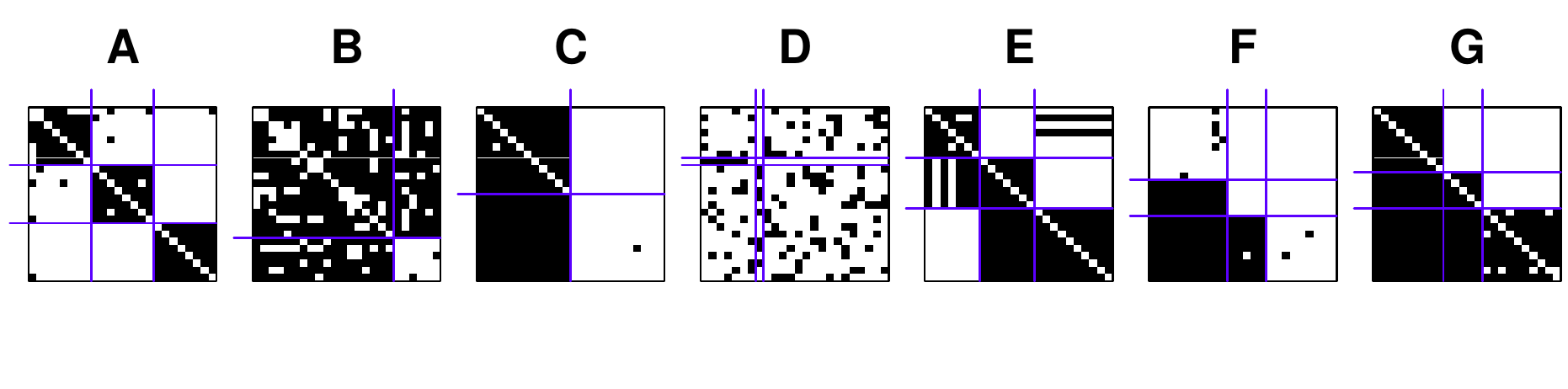}
\end{adjustwidth}
\end{figure}

Both algorithms perform well in the case of asymmetric core-periphery blockmodels (all generated networks have the expected blockmodel, see Fig~\ref{Fig10}). However, in the case of the symmetric core-periphery blockmodel, the MIVs are usually small, which is reflected by the insufficiently small periphery in the generated networks. It is also hard to generate a hierarchical blockmodel without complete blocks on the diagonal when considering only different triad types, regardless of the algorithm used to generate the networks. By adding paths of length three, the empirical networks produced have the expected blockmodel type with a very low level of errors (see Fig~\ref{Fig9}). 

One of the most important observations is that the number of different types of triads reflects the assumed global network structure, where it is often sufficient to consider only some of all possible types of triads. 
    
\section*{Conclusion}
    
This paper examined whether often used global structures, operationalised by blockmodels, can emerge as a consequence of local processes operationalised by different triad types. It was shown that for most of the studied blockmodels this indeed happened. The only exception is a hierarchical blockmodel without complete blocks on the diagonal, where additional local structures are needed to obtain a good fit with the assumed blockmodel structure. 

The main conclusion of the study is that a given global network structure can emerge due to the local mechanisms, regardless of the characteristics of the nodes. 

This was shown by generating networks that considered different triad types using two algorithms: the proposed deterministic Relocating Links (RL) algorithm and the Monte Carlo Markov Chain (MCMC) algorithm as implemented in the ergm R package \cite{hunter2008ergm}. The RL algorithm randomly selects a link and exchanges it with a randomly selected non-link. The change is accepted if the new network's local structure count is closer to the target count than in the previous network. With the MCMC algorithm, the same local structures were used as parameters in the ergm model.

To determine the target count for the RL algorithm and the parameter values for the MCMC algorithm, the count of different triad types in ideal networks (namely, those that perfectly comply with a certain blockmodel) had to be determined. This was achieved by considering the specific blockmodel type and corresponding group sizes. All types of triads were classified in the set of forbidden and set of allowed triad types for each blockmodel. Allowed triad types are those that are present in ideal networks and forbidden triad types are those that are not. The RL algorithm uses counts of a selected local structure in an ideal blockmodel while for the MCMC algorithm the parameter values were determined based on the classification into allowed and forbidden triad types.

Both algorithms performed very well; the exception is the hierarchical model without complete blocks on the diagonal (an additional parameter must be added) and to a smaller extent the symmetrical core-periphery model. On average, the RL algorithm performed slightly better. 

This paper also explored whether one can reduce the required local structure information by using only allowed or only forbidden triad types. Using only forbidden types of triads is especially desirable for the RL algorithm as the count for this triad type is zero. In addition, the reduction of all these sets (all, allowed, forbidden) of triad types was studied based on their sensitivity to errors according to the blockmodel structure. Most of these reductions of sets of different triad types overall resulted in only a slightly worse fit and in some cases even in an improved performance. The only exception is when using only selected forbidden triad types, which did not generate the assumed blockmodel structure.

There are several ways this study could be extended. One would be to move from considering the local structures (e.g. types of triads) to explicitly define other types of rules for creating and dissolving links in the network (local mechanisms). Another extension of this study could be to consider the mechanisms that lead to change from one particular blockmodel to another particular blockmodel.
    
\bibliographystyle{plos2015}
\bibliography{mybibfile}

\newpage
\begin{center}
\vspace*{\fill}
{\Huge SUPPLEMENTARY MATERIAL}
\vspace*{\fill}
\end{center}
\newpage

\vspace*{0.2in}
\begin{figure}[H]
\begin{adjustwidth}{-2.25in}{0in}
\caption*{{\bf S1 Fig. Some empirical generated networks using the RL algorithm by considering all triad types.} By rows: (A) cohesive; (B) symmetric core-periphery; (C) asymmetric core-periphery; (D) hierarchical without complete blocks on the diagonal; (E) hierarchical with complete blocks on the diagonal; (F) transitivity without complete blocks on the diagonal; (G) transitivity with complete blocks on the diagonal.}
    \label{emp}
    \includegraphics[width=1.3\textwidth]{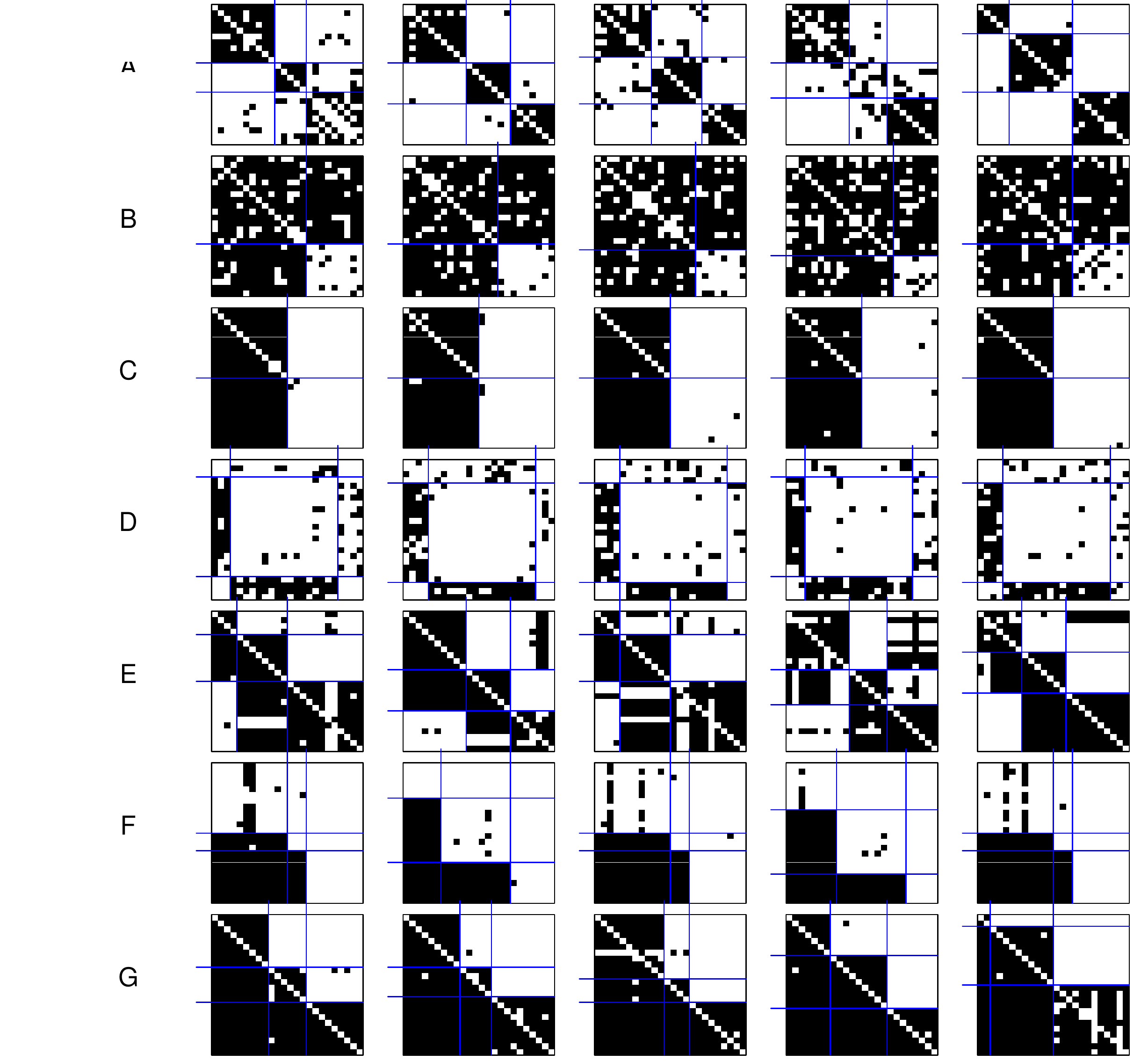}
\end{adjustwidth}
\end{figure}

\begin{figure}[H]
\begin{adjustwidth}{-2.25in}{0in}
\caption*{{\bf S2 Fig. Some empirical generated networks using the RL algorithm by considering allowed triad types.} By rows: (A) cohesive; (B) symmetric core-periphery; (C) asymmetric core-periphery; (D) hierarchical without complete blocks on the diagonal; (E) hierarchical with complete blocks on the diagonal; (F) transitivity without complete blocks on the diagonal; (G) transitivity with complete blocks on the diagonal.}
    \label{emp}
    \includegraphics[width=1.3\textwidth]{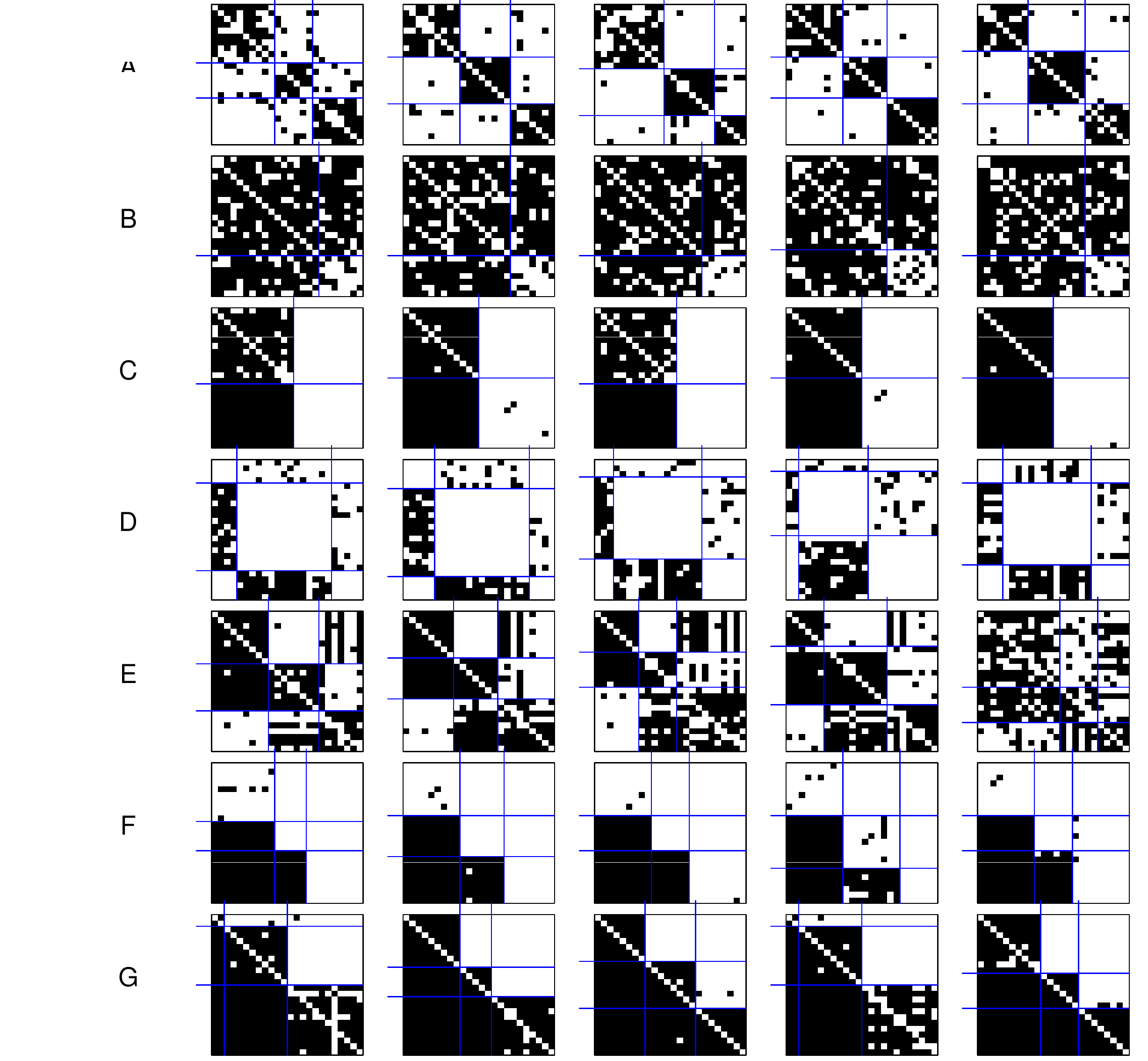}
\end{adjustwidth}
\end{figure}

\begin{figure}[H]
\begin{adjustwidth}{-2.25in}{0in}
\caption*{{\bf S3 Fig. Some empirical generated networks using the RL algorithm by considering forbidden triad types.} By rows: (A) cohesive; (B) symmetric core-periphery; (C) asymmetric core-periphery; (D) hierarchical without complete blocks on the diagonal; (E) hierarchical with complete blocks on the diagonal; (F) transitivity without complete blocks on the diagonal; (G) transitivity with complete blocks on the diagonal.}
    \label{emp}
    \includegraphics[width=1.3\textwidth]{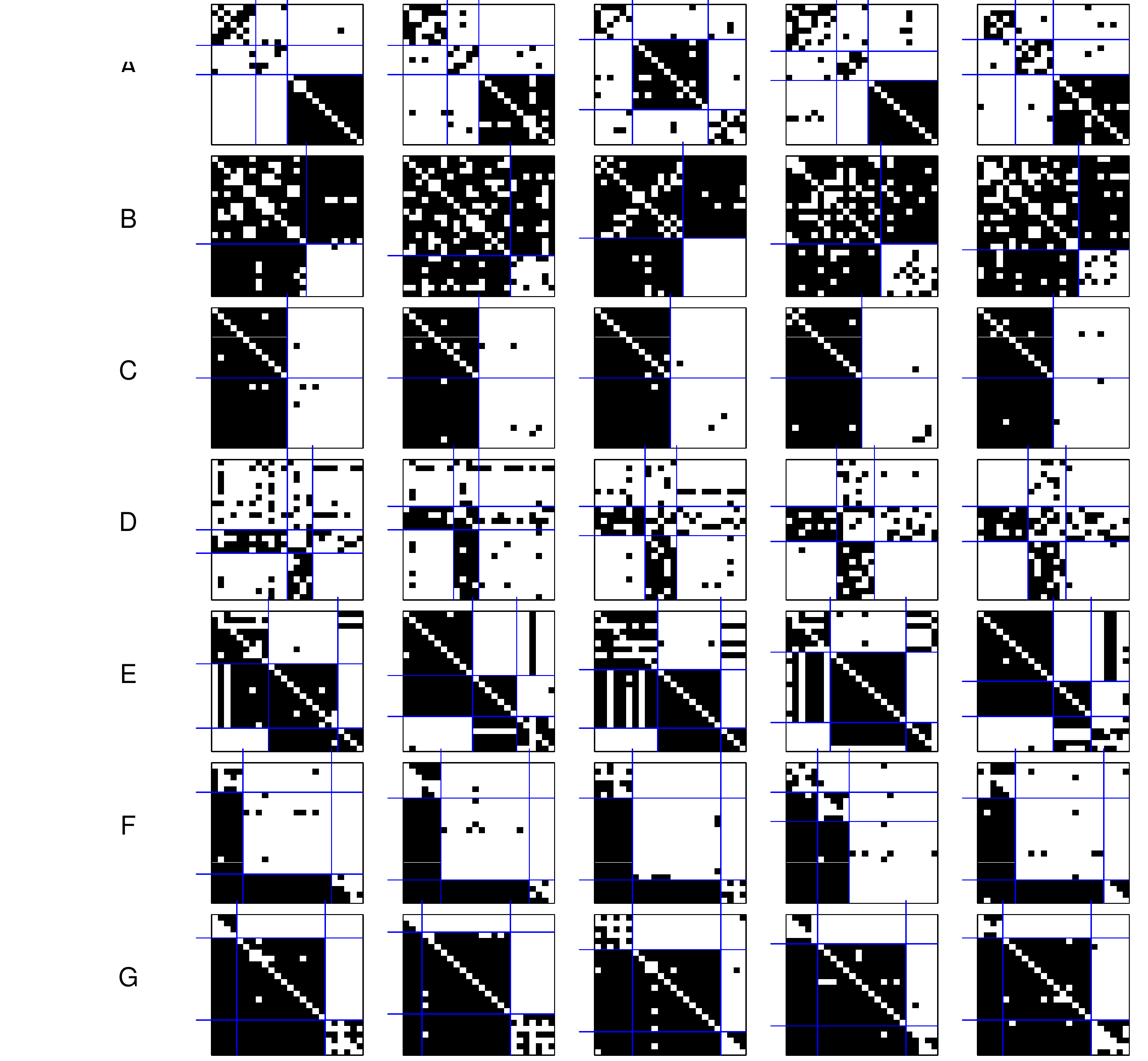}
\end{adjustwidth}
\end{figure}

\begin{figure}[H]
\begin{adjustwidth}{-2.25in}{0in}
\caption*{{\bf S4 Fig. Some empirical generated networks using the RL algorithm by considering all selected triad types.} By rows: (A) cohesive; (B) symmetric core-periphery; (C) asymmetric core-periphery; (D) hierarchical without complete blocks on the diagonal; (E) hierarchical with complete blocks on the diagonal; (F) transitivity without complete blocks on the diagonal; (G) transitivity with complete blocks on the diagonal.}
    \label{emp}
    \includegraphics[width=1.3\textwidth]{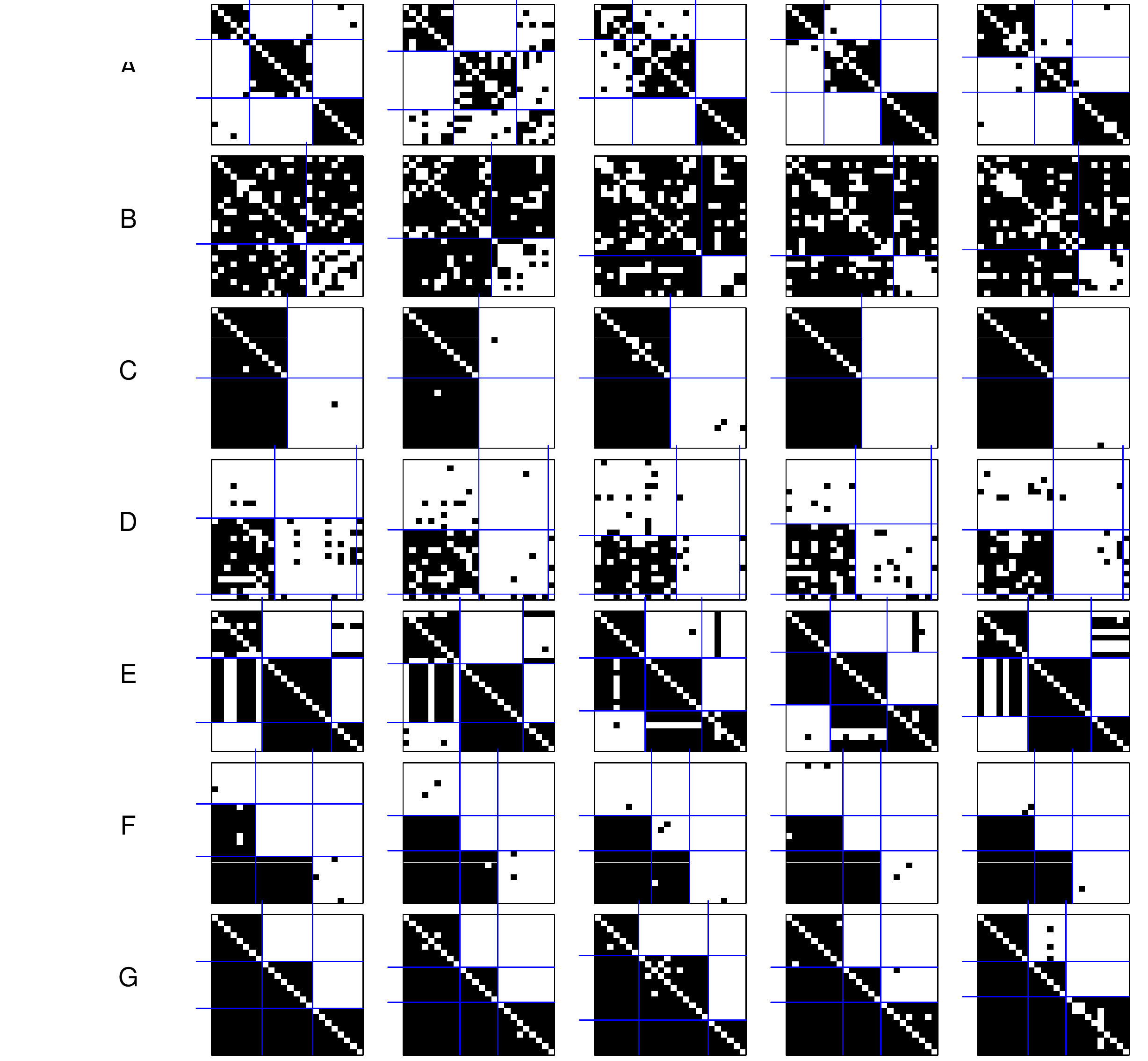}
\end{adjustwidth}
\end{figure}

\begin{figure}[H]
\begin{adjustwidth}{-2.25in}{0in}
\caption*{{\bf S5 Fig. Some empirical generated networks using the RL algorithm by considering selected allowed triad types.} By rows: (A) cohesive; (B) symmetric core-periphery; (C) asymmetric core-periphery; (D) hierarchical without complete blocks on the diagonal; (E) hierarchical with complete blocks on the diagonal; (F) transitivity without complete blocks on the diagonal; (G) transitivity with complete blocks on the diagonal.}
    \label{emp}
    \includegraphics[width=1.3\textwidth]{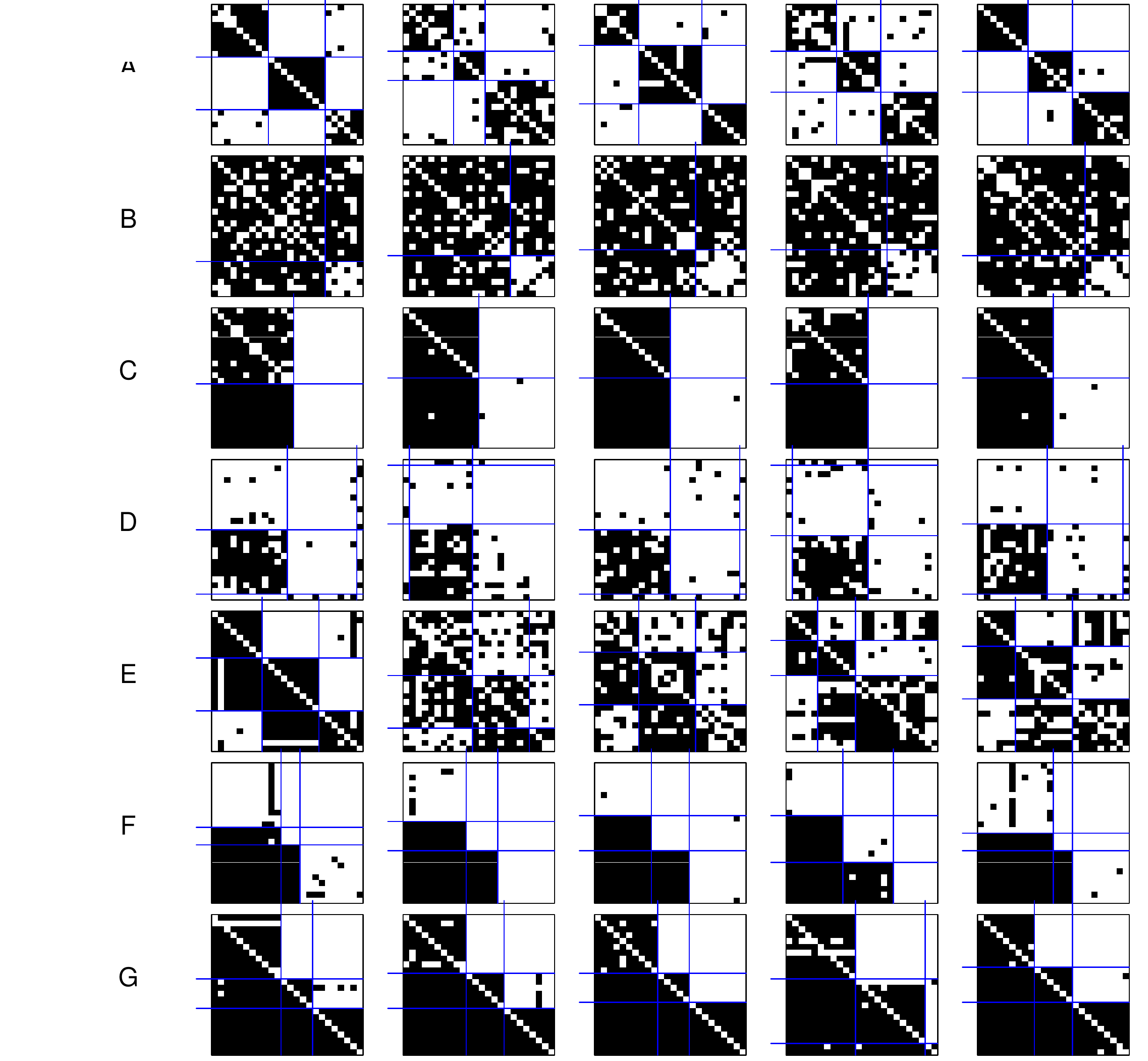}
\end{adjustwidth}
\end{figure}

\begin{figure}[H]
\begin{adjustwidth}{-2.25in}{0in}
\caption*{{\bf S6 Fig. Some empirical generated networks using the RL algorithm by considering selected forbidden triad types.} By rows: (A) cohesive; (B) symmetric core-periphery; (C) asymmetric core-periphery; (D) hierarchical without complete blocks on the diagonal; (E) hierarchical with complete blocks on the diagonal; (F) transitivity without complete blocks on the diagonal; (G) transitivity with complete blocks on the diagonal.}
    \label{emp}
    \includegraphics[width=1.3\textwidth]{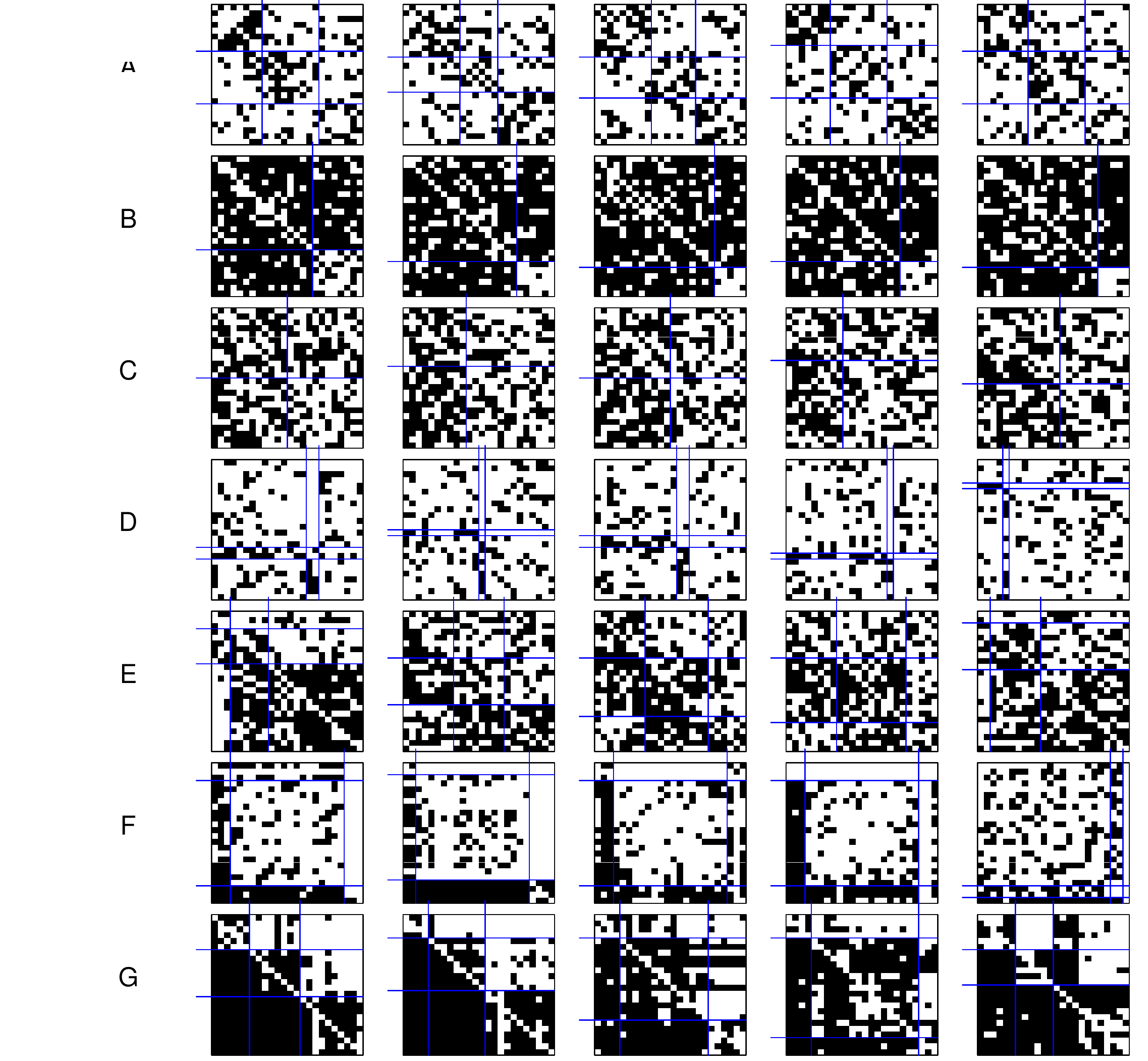}
\end{adjustwidth}
\end{figure}
       
\begin{figure}[H]
\begin{adjustwidth}{-2.25in}{0in}
\caption*{{\bf S7 Fig. Some empirical generated networks using the MCMC algorithm (fixed density) by considering all triad types.} By rows: (A) cohesive; (B) symmetric core-periphery; (C) asymmetric core-periphery; (D) hierarchical without complete blocks on the diagonal; (E) hierarchical with complete blocks on the diagonal; (F) transitivity without complete blocks on the diagonal; (G) transitivity with complete blocks on the diagonal.}
    \label{emp}
    \includegraphics[width=1.3\textwidth]{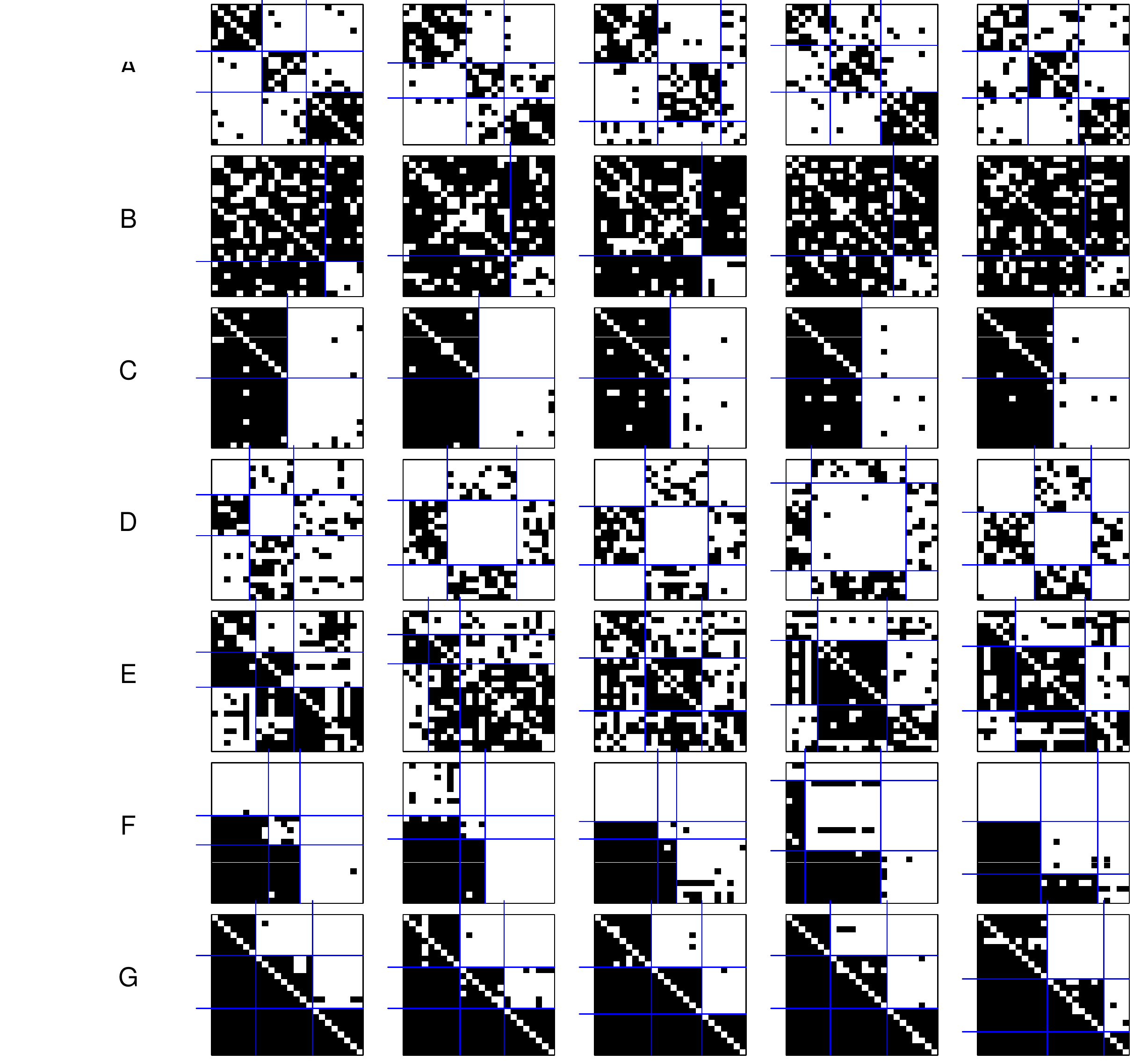}
\end{adjustwidth}
\end{figure}

\begin{figure}[H]
\begin{adjustwidth}{-2.25in}{0in}
\caption*{{\bf S8 Fig. Some empirical generated networks using the MCMC algorithm (fixed density) by considering allowed triad types.} By rows: (A) cohesive; (B) symmetric core-periphery; (C) asymmetric core-periphery; (D) hierarchical without complete blocks on the diagonal; (E) hierarchical with complete blocks on the diagonal; (F) transitivity without complete blocks on the diagonal; (G) transitivity with complete blocks on the diagonal.}
    \label{emp}
    \includegraphics[width=1.3\textwidth]{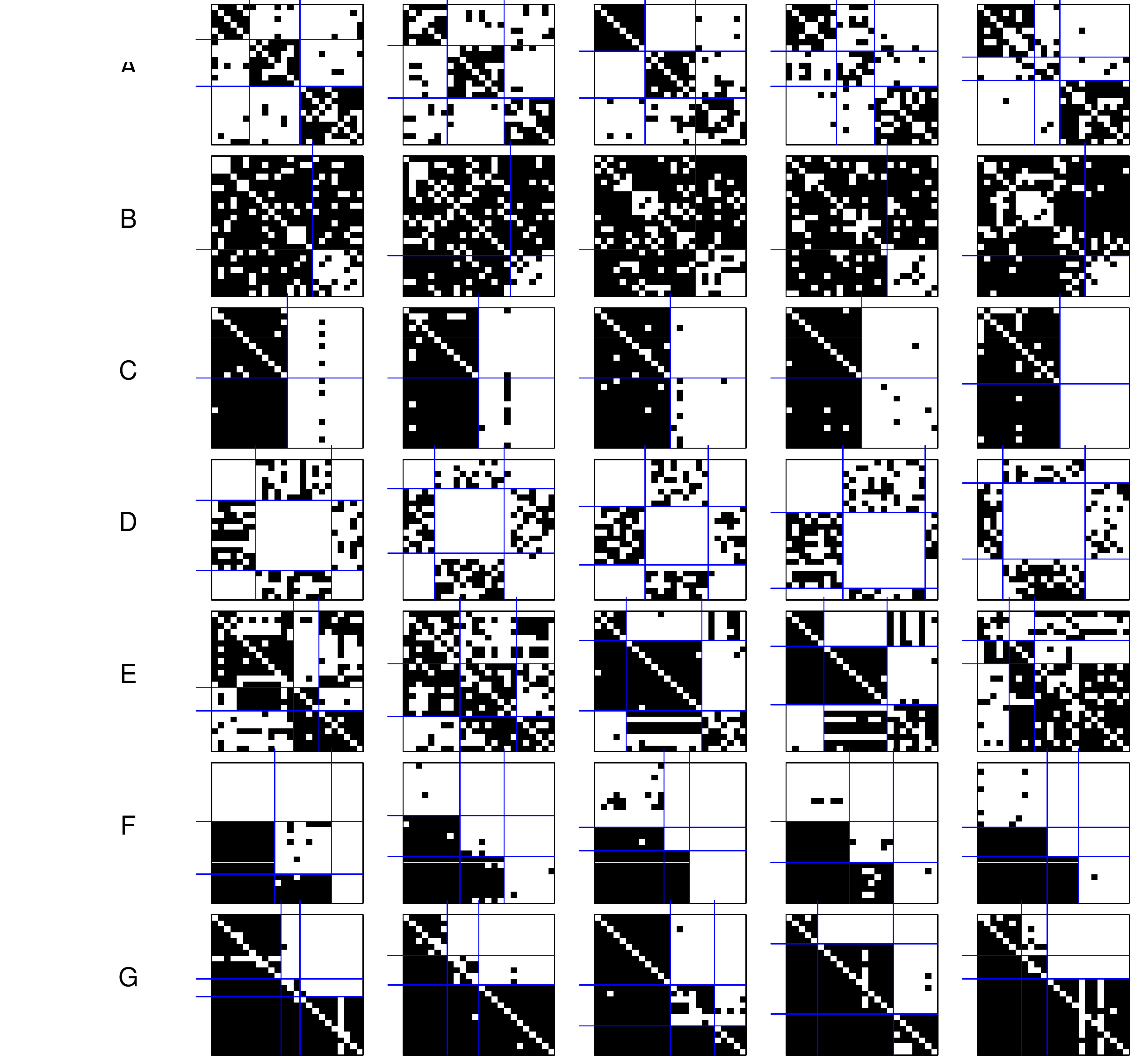}
\end{adjustwidth}
\end{figure}

\begin{figure}[H]
\begin{adjustwidth}{-2.25in}{0in}
\caption*{{\bf S9 Fig. Some empirical generated networks using the MCMC algorithm (fixed density) by considering forbidden triad types.} By rows: (A) cohesive; (B) symmetric core-periphery; (C) asymmetric core-periphery; (D) hierarchical without complete blocks on the diagonal; (E) hierarchical with complete blocks on the diagonal; (F) transitivity without complete blocks on the diagonal; (G) transitivity with complete blocks on the diagonal.}
    \label{emp}
    \includegraphics[width=1.3\textwidth]{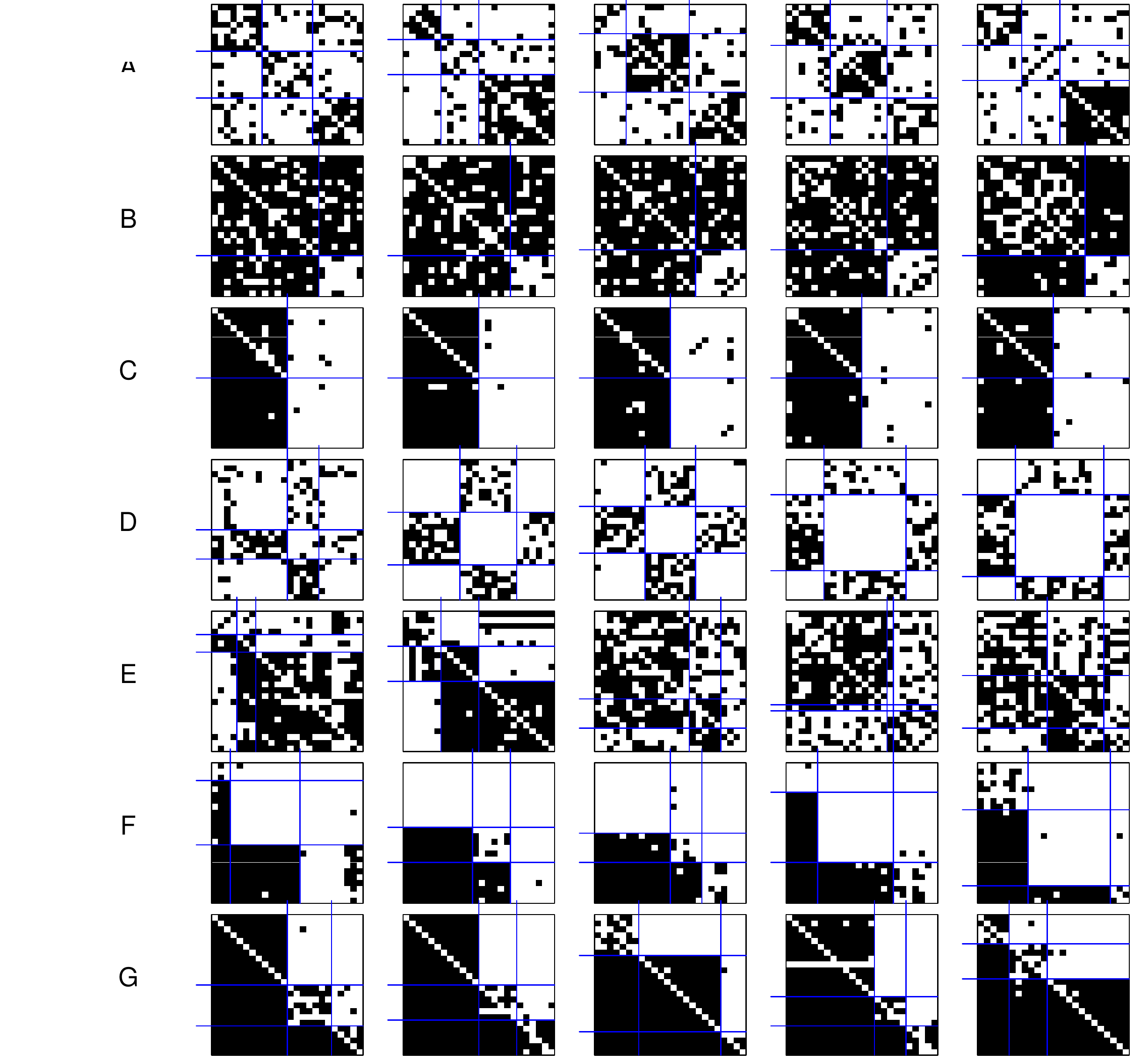}
\end{adjustwidth}
\end{figure}

\begin{figure}[H]
\begin{adjustwidth}{-2.25in}{0in}
\caption*{{\bf S10 Fig. Some empirical generated networks using the MCMC algorithm (fixed density) by considering all selected triad types.} By rows: (A) cohesive; (B) symmetric core-periphery; (C) asymmetric core-periphery; (D) hierarchical without complete blocks on the diagonal; (E) hierarchical with complete blocks on the diagonal; (F) transitivity without complete blocks on the diagonal; (G) transitivity with complete blocks on the diagonal.}
    \label{emp}
    \includegraphics[width=1.3\textwidth]{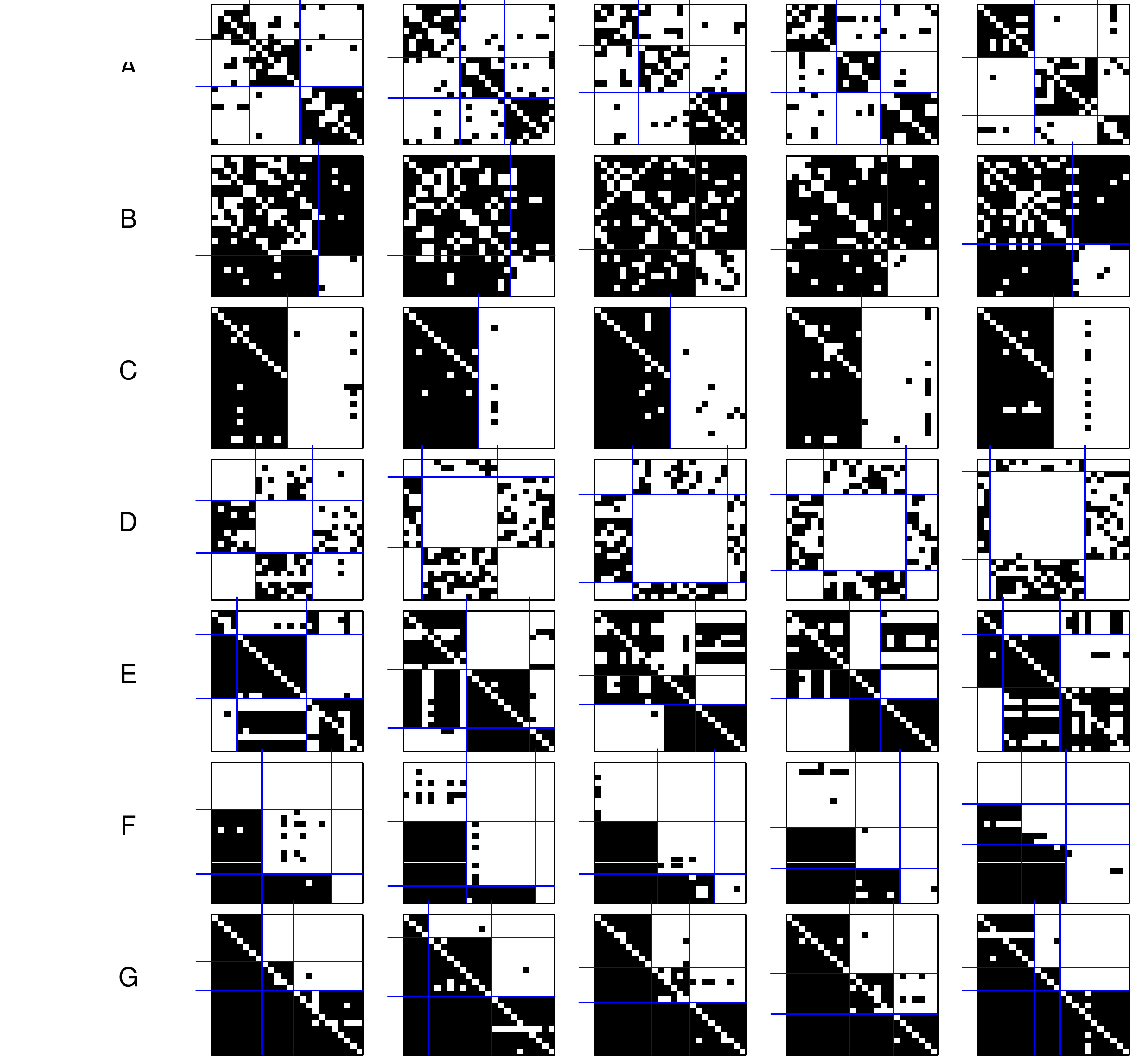}
\end{adjustwidth}
\end{figure}

\begin{figure}[H]
\begin{adjustwidth}{-2.25in}{0in}
\caption*{{\bf S11 Fig. Some empirical generated networks using the MCMC algorithm (fixed density) by considering selected allowed triad types.} By rows: (A) cohesive; (B) symmetric core-periphery; (C) asymmetric core-periphery; (D) hierarchical without complete blocks on the diagonal; (E) hierarchical with complete blocks on the diagonal; (F) transitivity without complete blocks on the diagonal; (G) transitivity with complete blocks on the diagonal.}
    \label{emp}
    \includegraphics[width=1.3\textwidth]{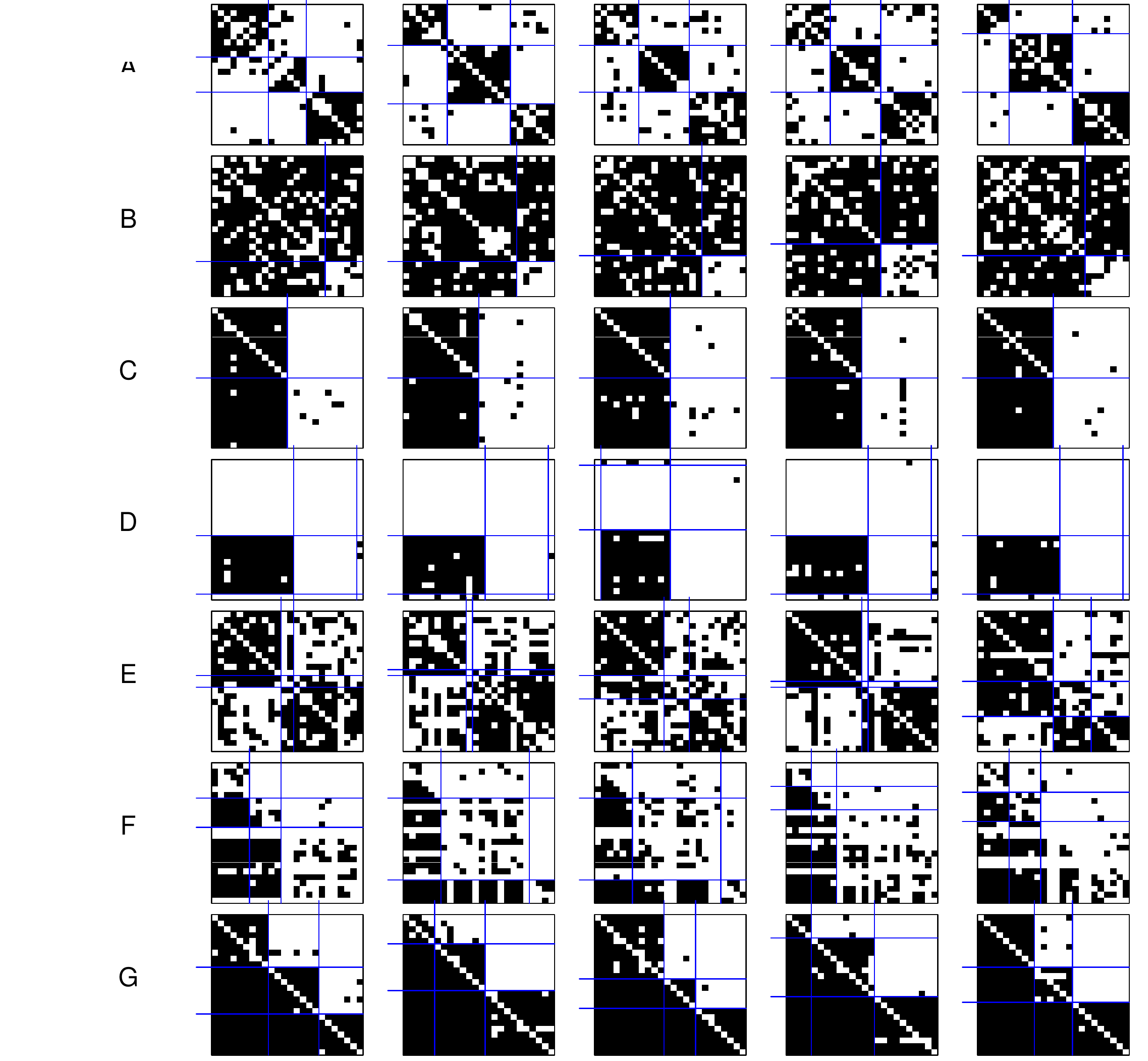}
\end{adjustwidth}
\end{figure}

\begin{figure}[H]
\begin{adjustwidth}{-2.25in}{0in}
\caption*{{\bf S12 Fig. Some empirical generated networks using the MCMC algorithm (fixed density) by considering selected forbidden triad types.} By rows: (A) cohesive; (B) symmetric core-periphery; (C) asymmetric core-periphery; (D) hierarchical without complete blocks on the diagonal; (E) hierarchical with complete blocks on the diagonal; (F) transitivity without complete blocks on the diagonal; (G) transitivity with complete blocks on the diagonal.}
    \label{emp}
    \includegraphics[width=1.3\textwidth]{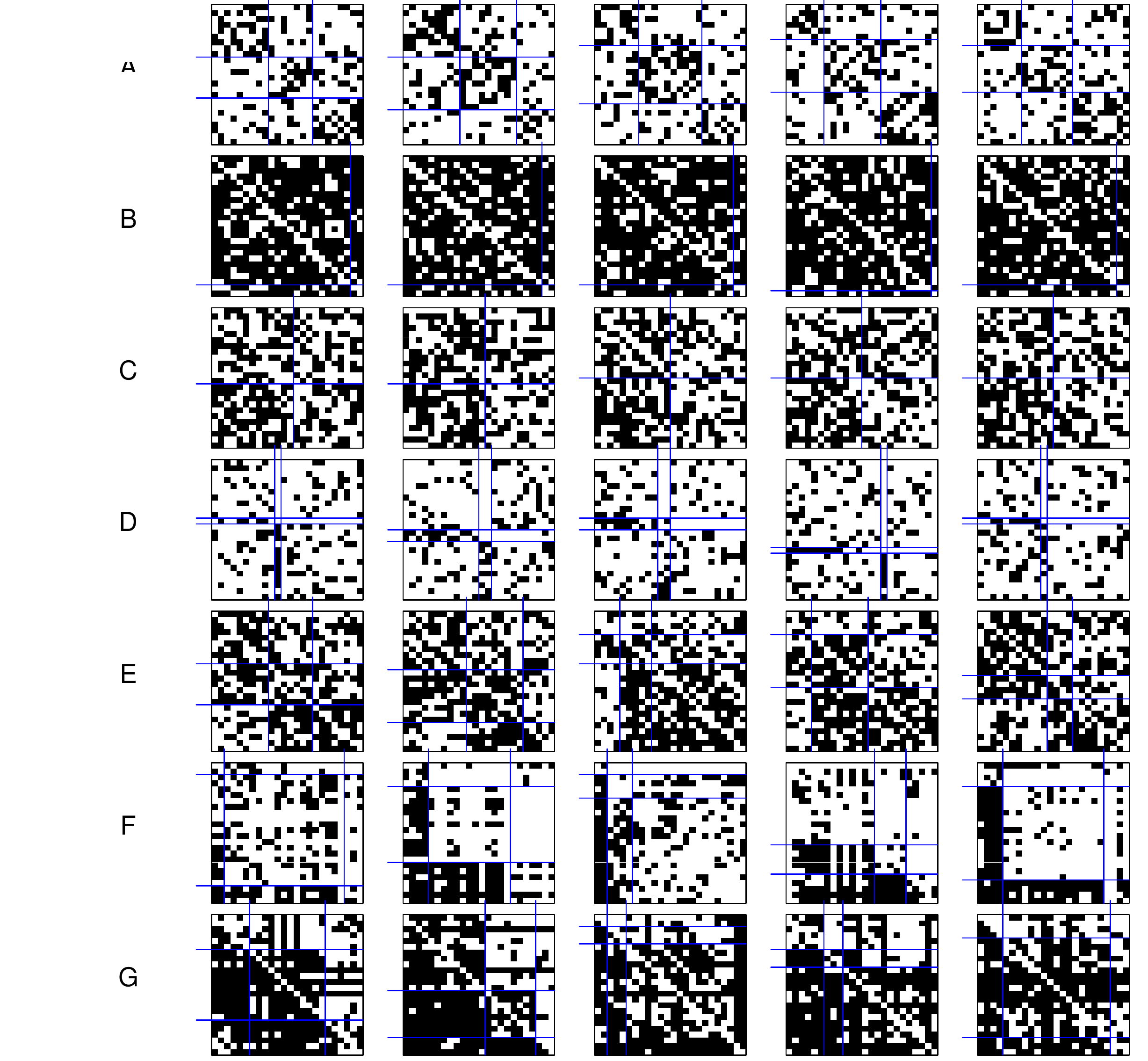}
\end{adjustwidth}
\end{figure}

\begin{figure}[H]
\begin{adjustwidth}{-2.25in}{0in}
\caption*{{\bf S13 Fig. Some empirical generated networks using the MCMC algorithm (variable density) by considering all triad types.} By rows: (A) cohesive; (B) symmetric core-periphery; (C) asymmetric core-periphery; (D) hierarchical without complete blocks on the diagonal; (E) hierarchical with complete blocks on the diagonal; (F) transitivity without complete blocks on the diagonal; (G) transitivity with complete blocks on the diagonal.}
    \label{emp}
    \includegraphics[width=1.3\textwidth]{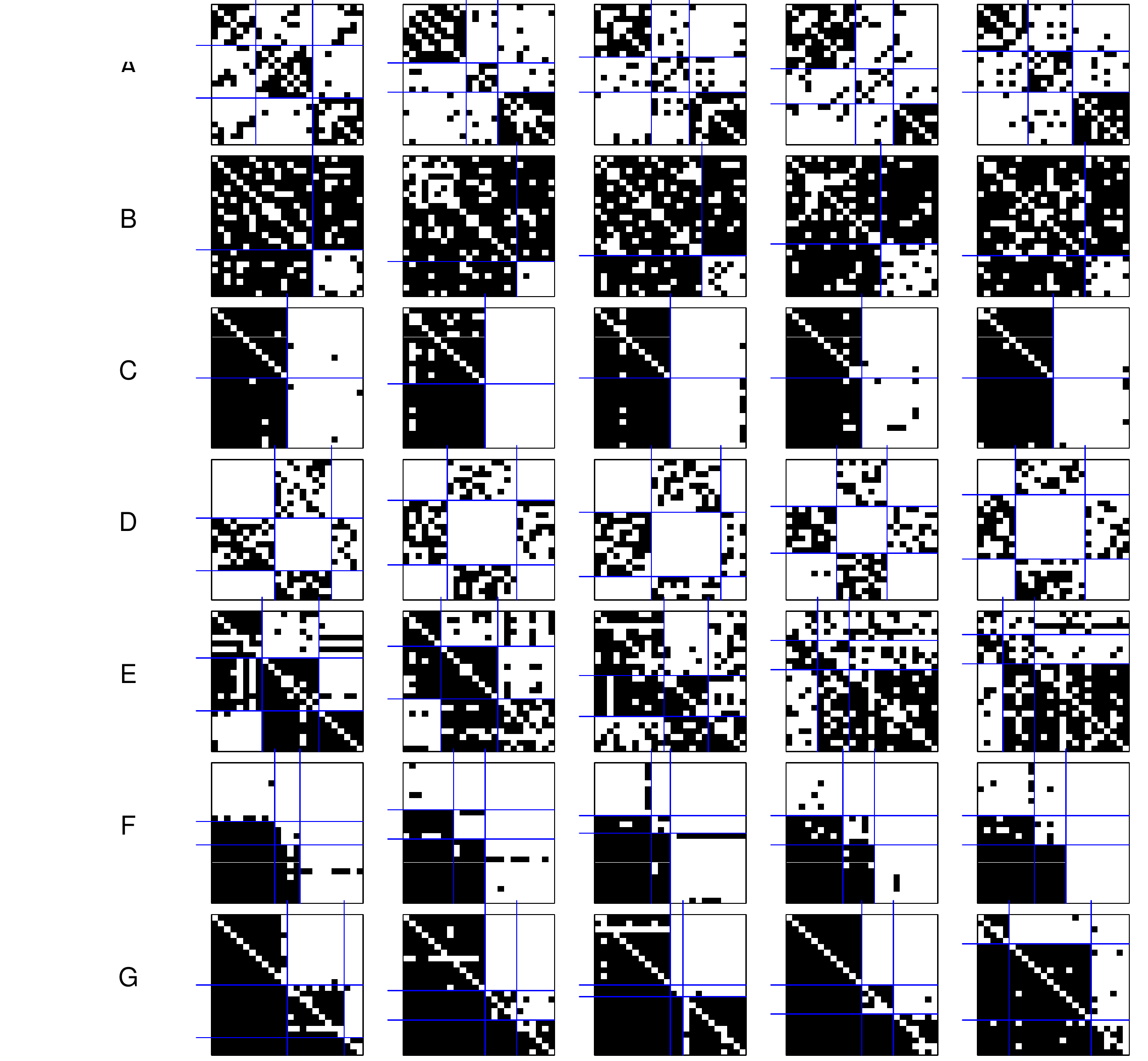}
\end{adjustwidth}
\end{figure}

\begin{figure}[H]
\begin{adjustwidth}{-2.25in}{0in}
\caption*{{\bf S14 Fig. Some empirical generated networks using the MCMC algorithm (variable density) by considering allowed triad types.} By rows: (A) cohesive; (B) symmetric core-periphery; (C) asymmetric core-periphery; (D) hierarchical without complete blocks on the diagonal; (E) hierarchical with complete blocks on the diagonal; (F) transitivity without complete blocks on the diagonal; (G) transitivity with complete blocks on the diagonal.}
    \label{emp}
    \includegraphics[width=1.3\textwidth]{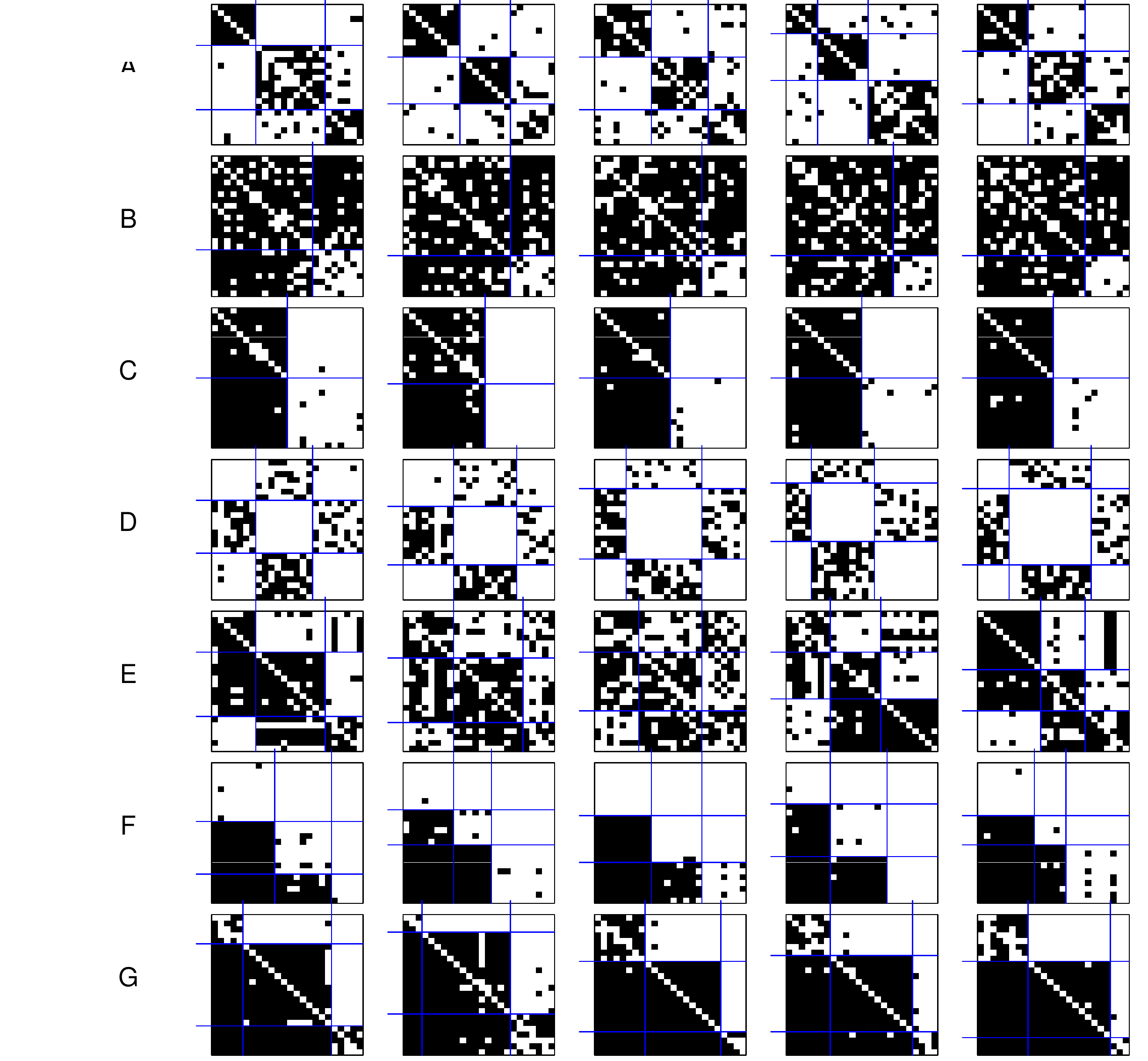}
\end{adjustwidth}
\end{figure}

\begin{figure}[H]
\begin{adjustwidth}{-2.25in}{0in}
\caption*{{\bf S15 Fig. Some empirical generated networks using the MCMC algorithm (variable density) by considering forbidden triad types.} By rows: (A) cohesive; (B) symmetric core-periphery; (C) asymmetric core-periphery; (D) hierarchical without complete blocks on the diagonal; (E) hierarchical with complete blocks on the diagonal; (F) transitivity without complete blocks on the diagonal; (G) transitivity with complete blocks on the diagonal.}
    \label{emp}
    \includegraphics[width=1.3\textwidth]{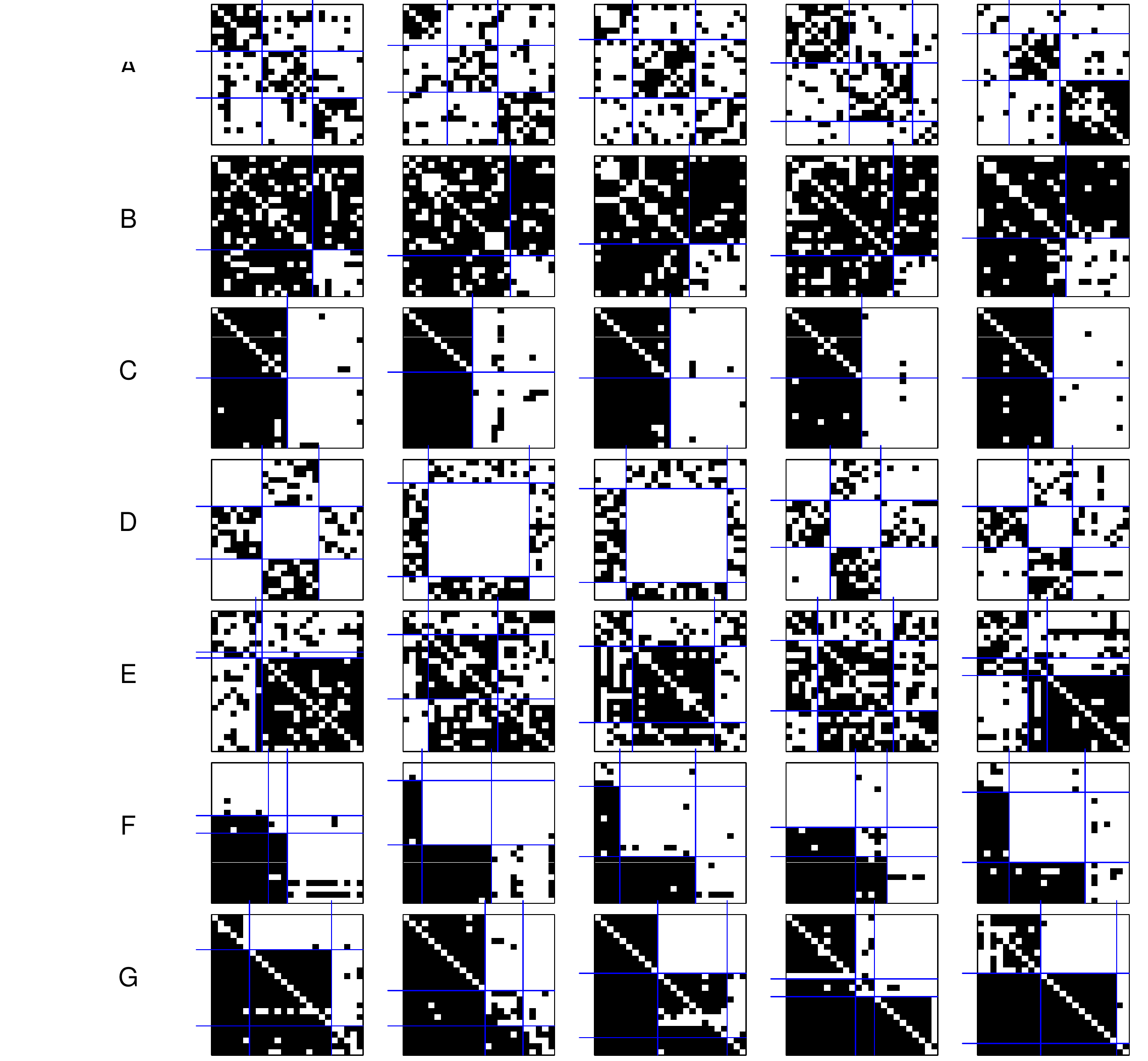}
\end{adjustwidth}
\end{figure}

\begin{figure}[H]
\begin{adjustwidth}{-2.25in}{0in}
\caption*{{\bf S16 Fig. Some empirical generated networks using the MCMC algorithm (variable density) by considering all selected triad types.} By rows: (A) cohesive; (B) symmetric core-periphery; (C) asymmetric core-periphery; (D) hierarchical without complete blocks on the diagonal; (E) hierarchical with complete blocks on the diagonal; (F) transitivity without complete blocks on the diagonal; (G) transitivity with complete blocks on the diagonal.}
    \label{emp}
    \includegraphics[width=1.3\textwidth]{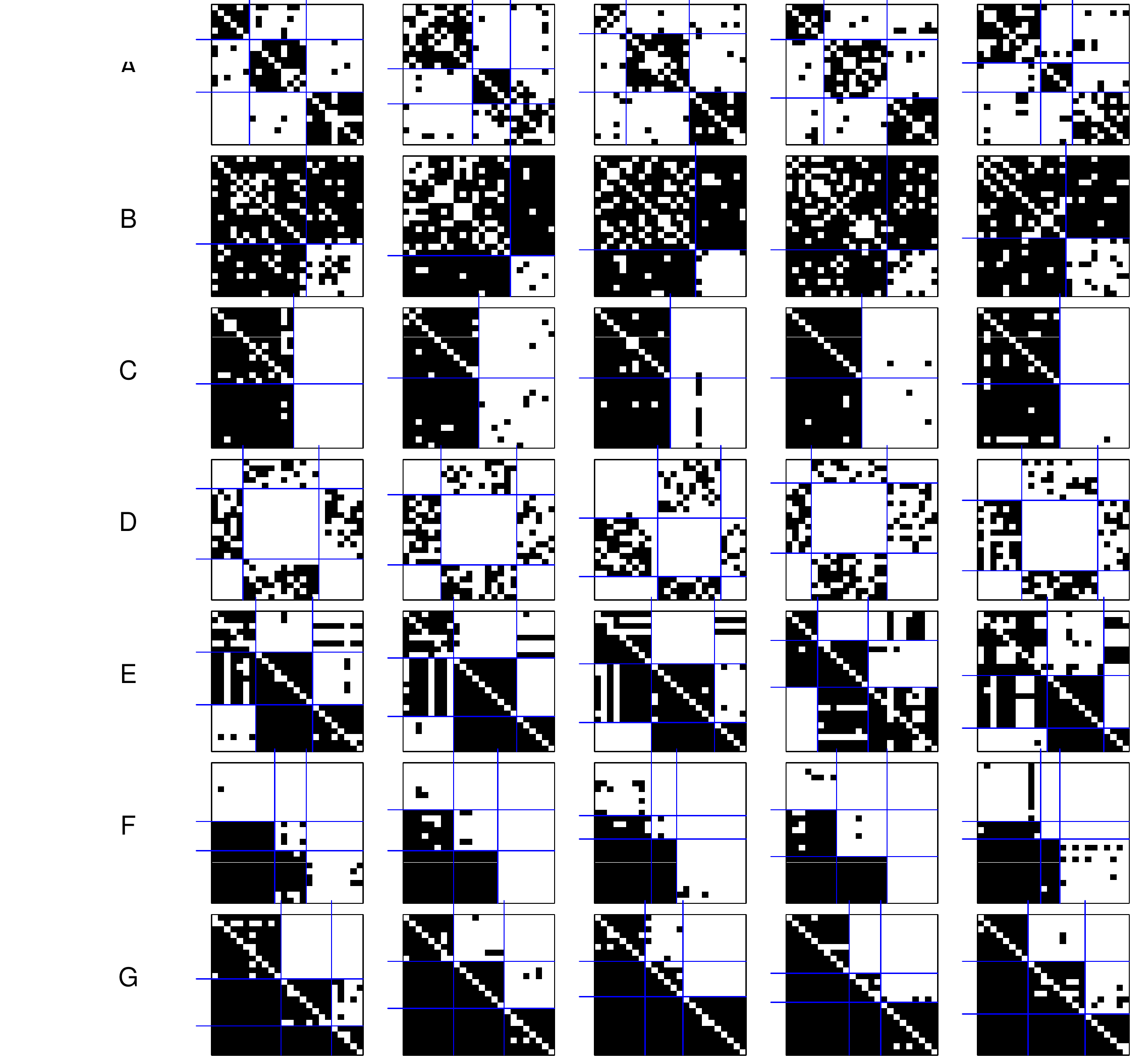}
\end{adjustwidth}
\end{figure}

\begin{figure}[H]
\begin{adjustwidth}{-2.25in}{0in}
\caption*{{\bf S17 Fig. Some empirical generated networks using the MCMC algorithm (variable density) by considering selected allowed triad types.} By rows: (A) cohesive; (B) symmetric core-periphery; (C) asymmetric core-periphery; (D) hierarchical without complete blocks on the diagonal; (E) hierarchical with complete blocks on the diagonal; (F) transitivity without complete blocks on the diagonal; (G) transitivity with complete blocks on the diagonal.}
    \label{emp}
    \includegraphics[width=1.3\textwidth]{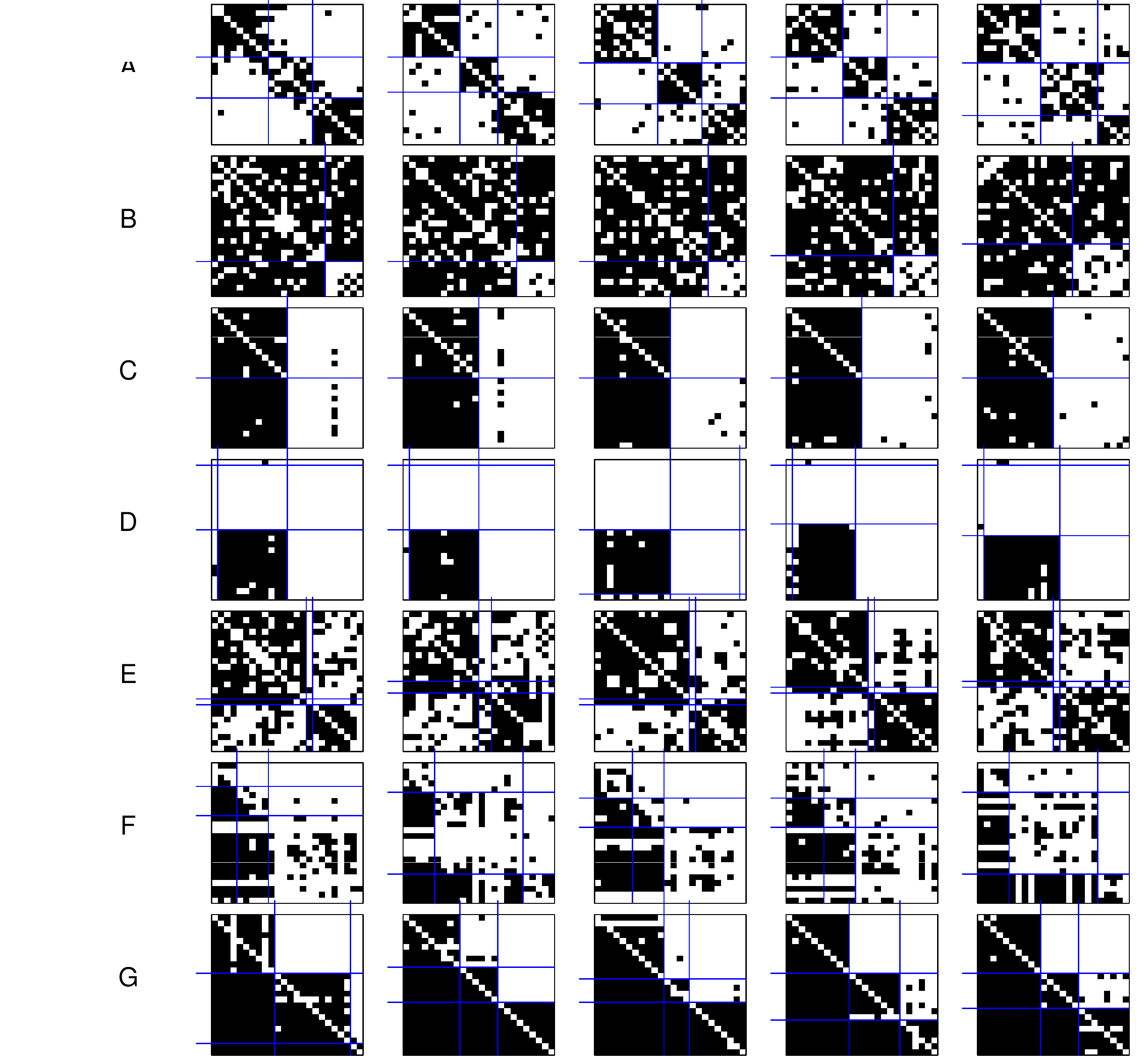}
\end{adjustwidth}
\end{figure}

\begin{figure}[H]
\begin{adjustwidth}{-2.25in}{0in}
\caption*{{\bf S18 Fig. Some empirical generated networks using the MCMC algorithm (variable density) by considering selected forbidden triad types.} By rows: (A) cohesive; (B) symmetric core-periphery; (C) asymmetric core-periphery; (D) hierarchical without complete blocks on the diagonal; (E) hierarchical with complete blocks on the diagonal; (F) transitivity without complete blocks on the diagonal; (G) transitivity with complete blocks on the diagonal.}
    \label{emp}
    \includegraphics[width=1.3\textwidth]{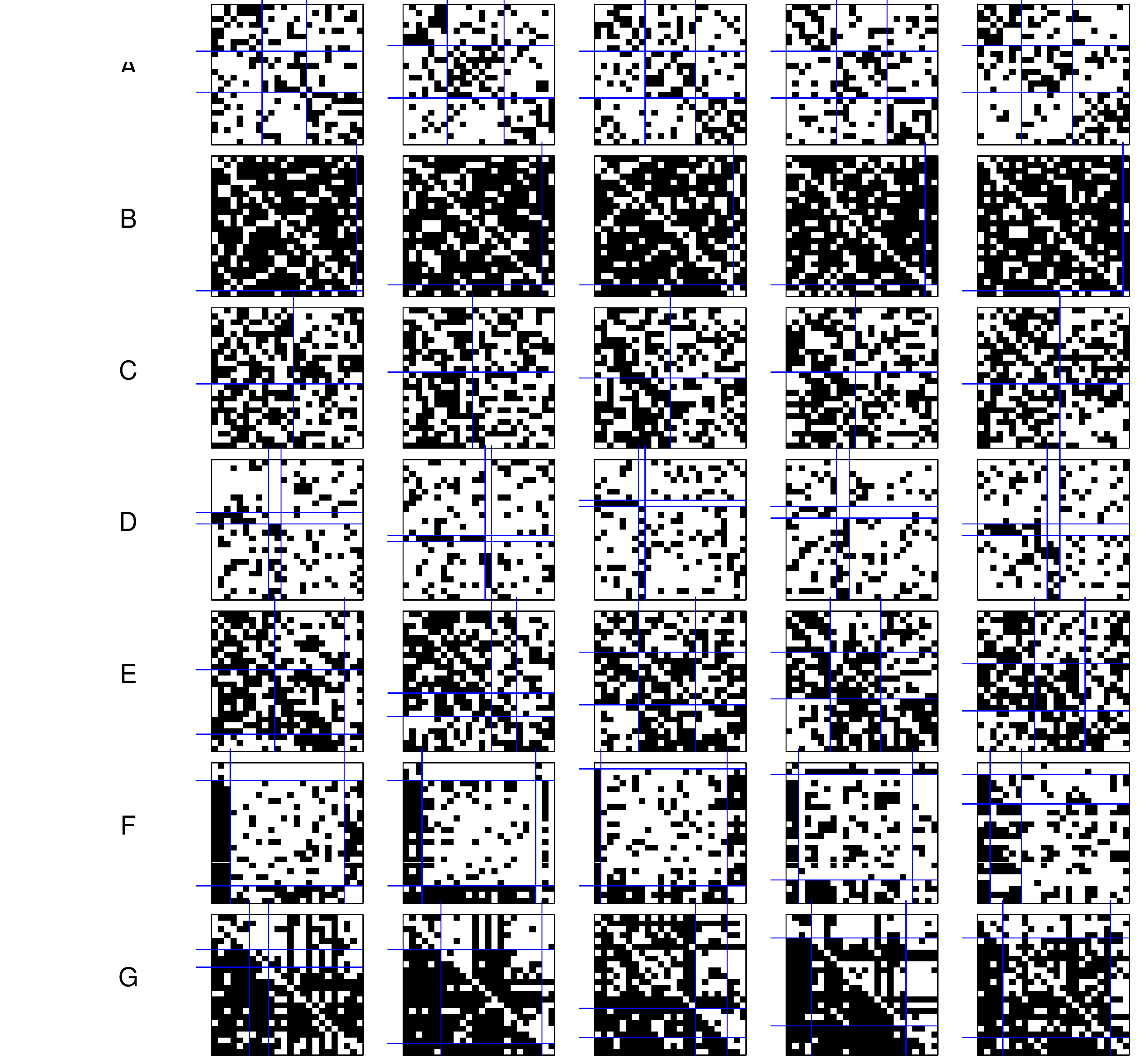}
\end{adjustwidth}
\end{figure}

    \begin{figure}[H]
\begin{adjustwidth}{-2.25in}{0in}
  \caption*{{\bf S19 Fig. The distribution of the criterion function for the generated networks and for randomised networks with the mean improvement value obtained by the RL algorithm.} R denotes the distribution of $P_R$ and M denotes the distribution of $P_m$. By rows: (A) cohesive; (B) symmetric core-periphery; (C) asymmetric core-periphery; (D) hierarchical without complete blocks on the diagonal; (E) hierarchical with complete blocks on the diagonal; (F) transitivity without complete blocks on the diagonal; (G) transitivity with complete blocks on the diagonal. By columns: (H) all; (I) all allowed; (J) all forbidden; (K) selected; (L) selected allowed; (M) selected forbidden types of triads.}
  \label{rlcr}
  \includegraphics[width=1.5\textwidth, angle = 270]{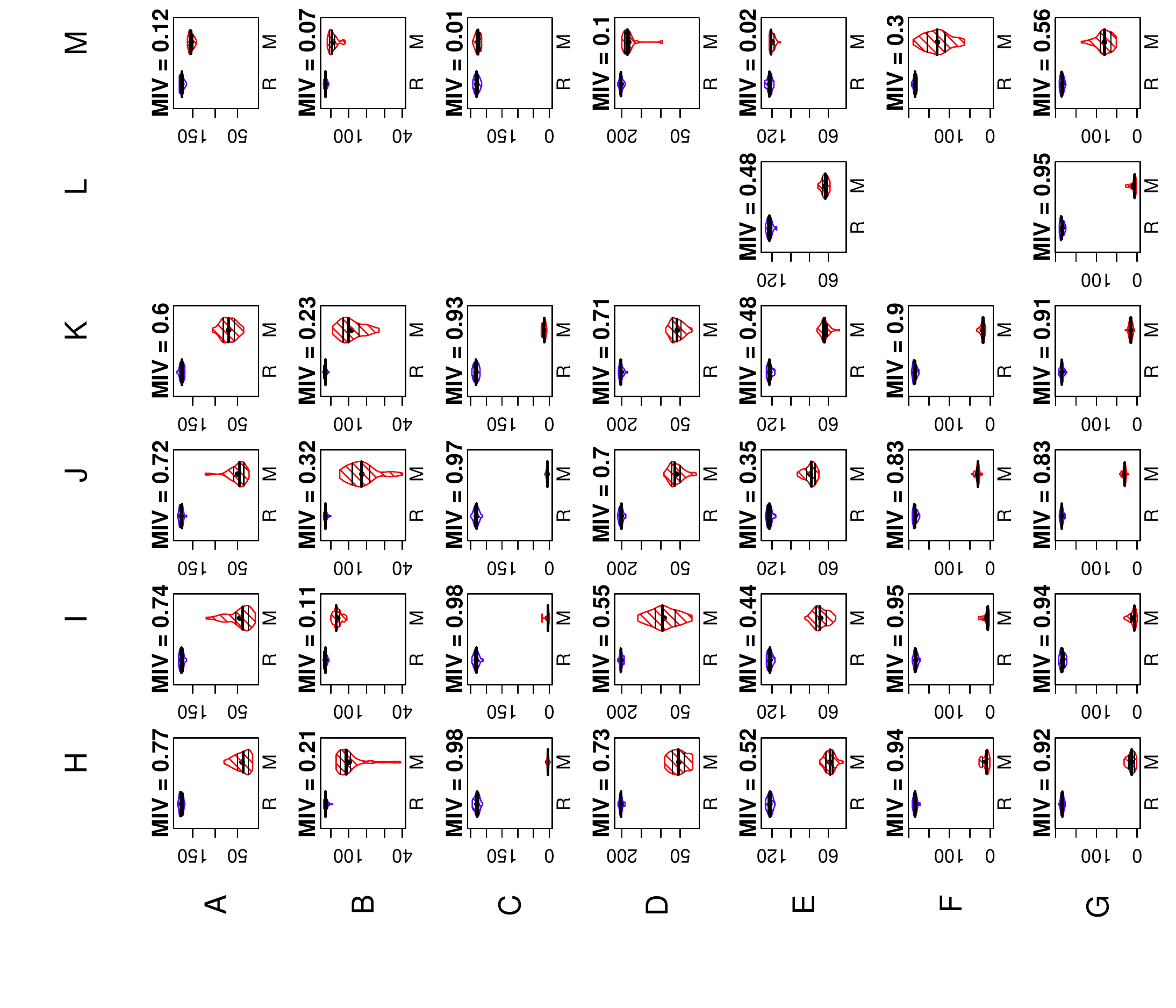}
\end{adjustwidth}
\end{figure}

\begin{figure}[H]
\begin{adjustwidth}{-2.25in}{0in}
  \caption*{{\bf S20 Fig. The distribution of the criterion function for the generated networks and for randomised networks with the mean improvement value obtained by the MCMC algorithm with fixed density.} R denotes the distribution of $P_R$ and M denotes the distribution of $P_m$. R denotes the distribution of $P_R$ and M denotes the distribution of $P_m$. By rows: (A) cohesive; (B) symmetric core-periphery; (C) asymmetric core-periphery; (D) hierarchical without complete blocks on the diagonal; (E) hierarchical with complete blocks on the diagonal; (F) transitivity without complete blocks on the diagonal; (G) transitivity with complete blocks on the diagonal. By columns: (H) all; (I) all allowed; (J) all forbidden; (K) selected; (L) selected allowed; (M) selected forbidden types of triads.}
  \label{mcmccrfg}
  \includegraphics[width=1.5\textwidth, angle = 270]{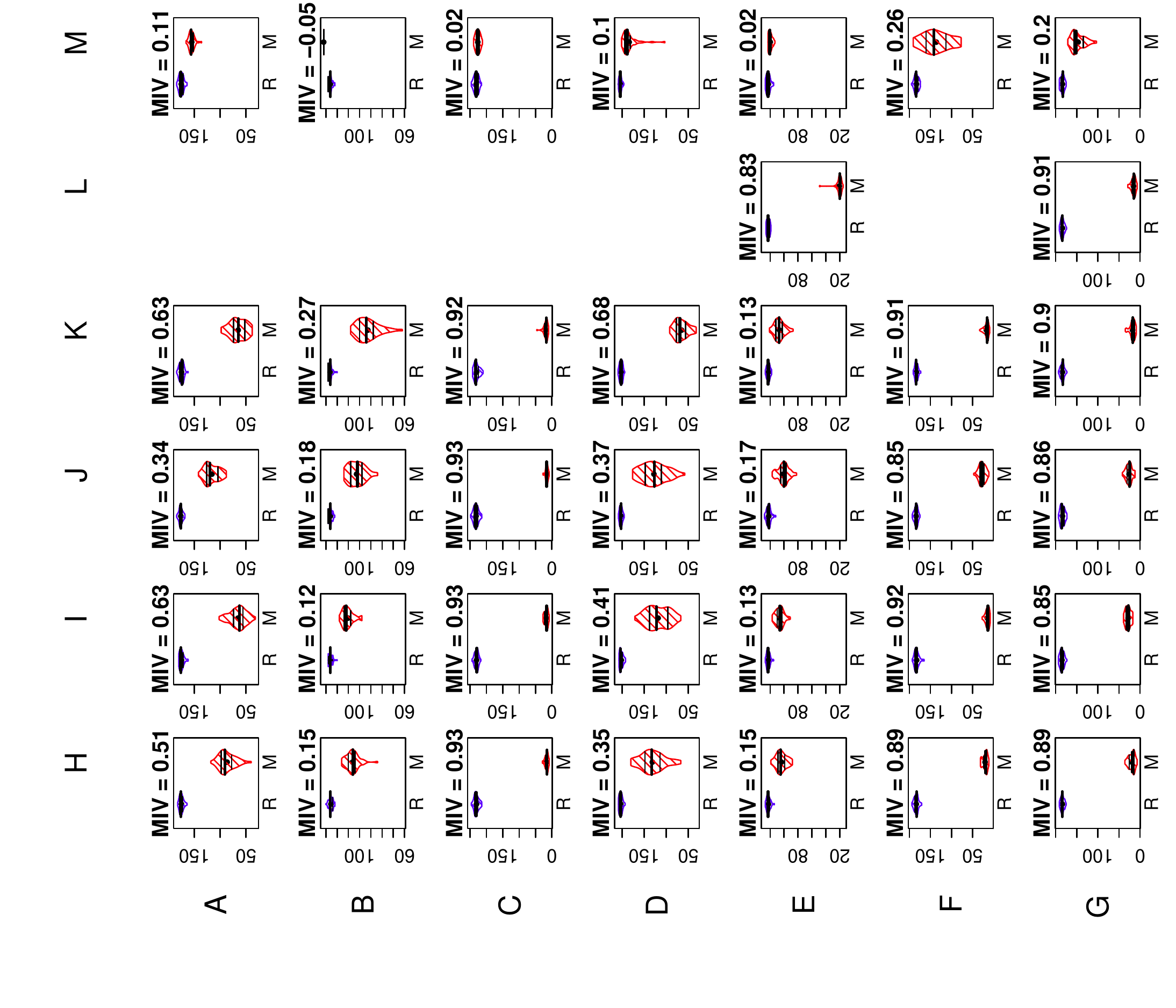}
\end{adjustwidth}
\end{figure}

\begin{figure}[H]
\begin{adjustwidth}{-2.25in}{0in}
  \caption*{{\bf S21 Fig. The distribution of the criterion function for the generated networks and for randomised networks with the mean improvement value obtained by the MCMC algorithm with non-fixed density} R denotes the distribution of $P_R$ and M denotes the distribution of $P_m$. R denotes the distribution of $P_R$ and M denotes the distribution of $P_m$. By rows: (A) cohesive; (B) symmetric core-periphery; (C) asymmetric core-periphery; (D) hierarchical without complete blocks on the diagonal; (E) hierarchical with complete blocks on the diagonal; (F) transitivity without complete blocks on the diagonal; (G) transitivity with complete blocks on the diagonal. By columns: (H) all; (I) all allowed; (J) all forbidden; (K) selected; (L) selected allowed; (M) selected forbidden types of triads.}
  \label{mcmccr}
  \includegraphics[width=1.5\textwidth, angle = 270]{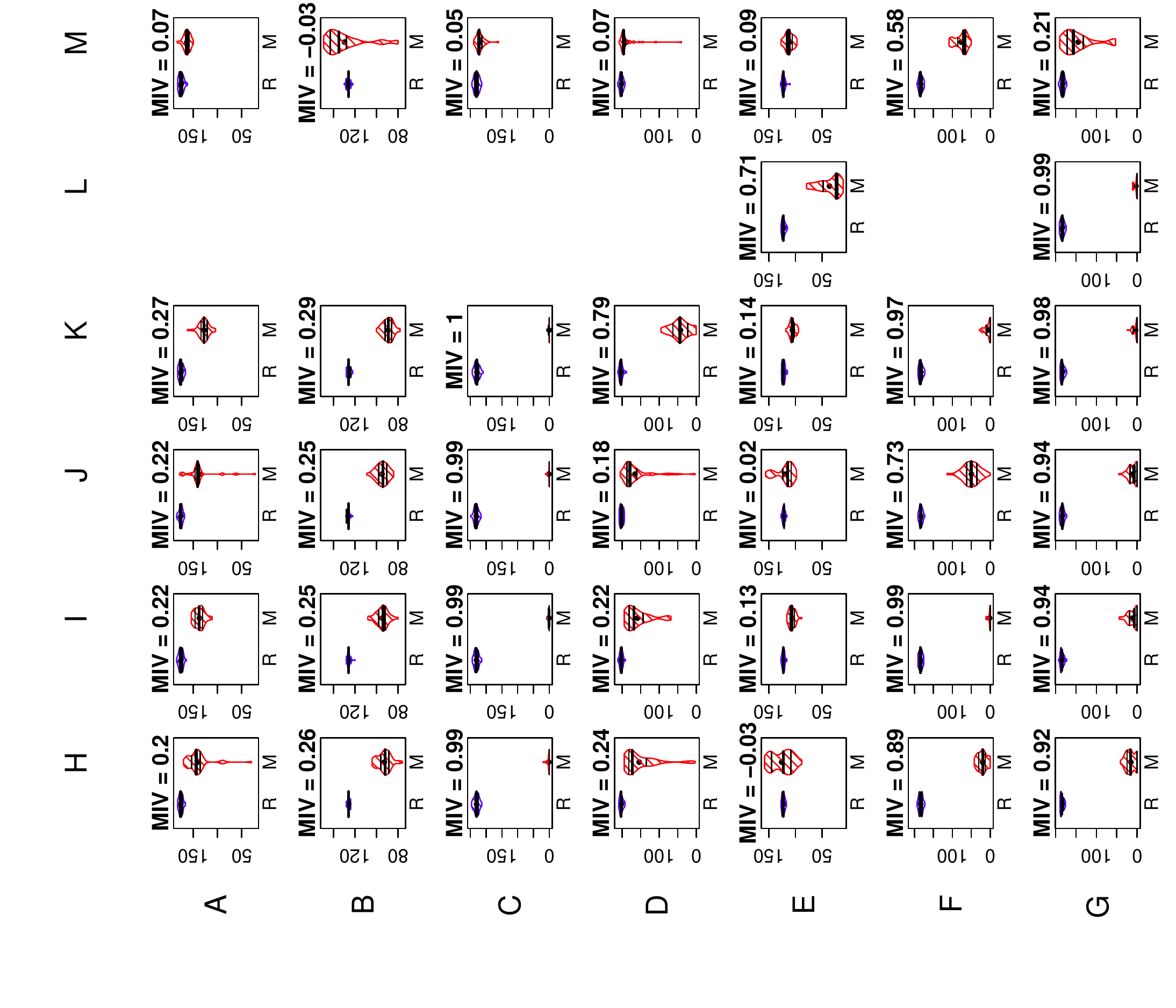}
\end{adjustwidth}
\end{figure}

\section*{S22 Generating totally randomised networks and networks with a given level of errors}

Totally randomised networks can be generated based on an ideal network: $k$ links in complete blocks are randomly chosen and replaced by non-links. At the same time, $k$ non-links are randomly chosen and replaced with links. In other words, the number of links is relocated in such a way that the overall density in complete blocks and overall density in null blocks are equal (i.e., all expected densities of all the blocks in the blockmodel are equal). The number of relocated links $k$ is calculated as 

\begin{eqnarray}
\label{k}
k=m - \frac{m^2}{n^2-n}
\end{eqnarray}

\noindent  where $m$ is the number of links and $n$ is the number of units in a selected type of blockmodel.

Instead of totally randomised networks, blockmodels with a certain level of errors can be analysed. In such case, when a network with a given blockmodel structure has to be generated with a certain level of errors, the number of relocated links is calculated as 

\begin{eqnarray}
\label{k2}
k={ m - \Big (}\frac{m^2}{n^2-n}{\Big )} * LE
\end{eqnarray}

\noindent where the level of errors ($LE$) can take a value on the interval $[0, 1]$ (0 stands for an ideal network and 1 for a random network). 

\section*{S23 Selected allowed and forbidden triad types}

The sets of allowed and forbidden triad types can be further reduced to selected allowed and selected forbidden triad types. This means that not all possible or all allowed or all forbidden types of triads are considered. One could choose the most appropriate terms based on observations of the A-measure. To obtain a better triad selection, it is beneficial to observe the A-values for networks with different levels of errors. 

\begin{figure}[H]
\begin{adjustwidth}{-2.25in}{0in}
\caption{{\bf A-measure values ($y$-axis) for different levels of errors ($x$-axis) and different types of blockmodels} (by rows: A = cohesive; B = symmetric core-periphery; C = asymmetric core-periphery; D = hierarchical without complete blocks on the diagonal; E = hierarchical with complete blocks on the diagonal; F = transitivity without complete blocks on the diagonal; G = transitivity with complete blocks on the diagonal)} 
\label{slika}
\centering
    \label{emp}
    \includegraphics[width=1.6\textwidth, angle = 90]{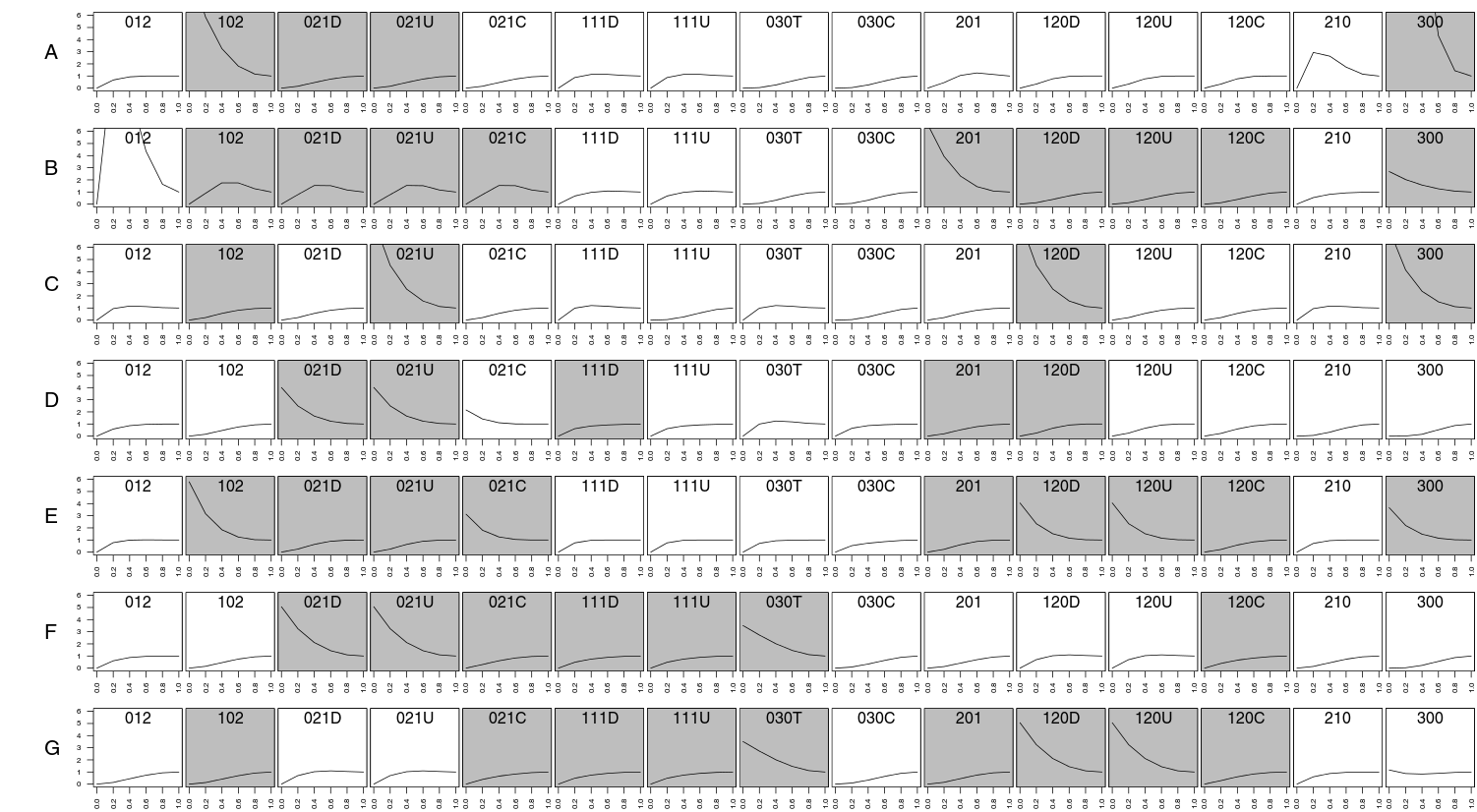}
\end{adjustwidth}
\end{figure}

Here, the most common (and uncommon) triads for each blockmodel type can be recognised by their sensitivity to different levels of errors. The idea is as follows: the most important triads are those with the highest absolute A-measure values for all levels of errors and with as close to a linear trend as possible through all levels of errors, indicating that a certain triad is not greatly affected by the level of errors (see Fig~\ref{slika} for the visualised relationship between the level of errors and A-measure for different blockmodel types and triad types). 

For some triads, the A-measure values are nearly constant for all levels of errors. Such a triad is triad type 300 in the case of a transitivity blockmodel with complete blocks on the diagonal. The value of the A-measure for these triad types is not associated with the level of errors. On the other hand, for many triad types a sharp change in the A-measure value at a certain level of errors is common. For example, in the case of a hierarchical blockmodel without complete blocks on the diagonal, the value of the A-measure for triads of types 012, 111D, 111U, 030T, 030C and 210 is zero in the case of an ideal network while it approaches 1 at very low levels of errors (i.e. between 0.2 and 0.4) and then remains constant. Some values first increase very fast at low levels of errors and then decrease at higher levels of errors. One example is the number of compete subgraphs of size three with one missing link in a cohesive blockmodel. 

The values of the A-measure for some types of triads are increasing or decreasing nearly linearly with the level of errors. These types of triads can be seen as triads that should be considered when generating networks with a given blockmodel structure. However, these types of triads can be further differentiated. For example, there are many types of triads with similar A-measure values for different levels of errors within some types of blockmodels. This could indicate that certain types of triads are defined similarly and are therefore not needed when generating networks with a given blockmodel structure. 

Some types of triads which are strongly influenced by the level of errors at low levels and less influenced by the level of errors at high level of errors (and vice versa) could also be chosen. In this case, it may happen that one should choose different types of triads for networks with higher and for networks with lower levels of errors.

\end{document}